\documentclass[journal]{IEEEtran}

\ifCLASSINFOpdf

\else

\fi

\usepackage{amsmath,graphicx}
\usepackage{amsmath,graphicx}
\usepackage{amsmath,amsfonts,amssymb,amsthm}
\usepackage{epstopdf}

\usepackage{url}

\DeclareMathOperator{\hist}{hist}
\DeclareMathOperator*{\argmin}{arg\,min}

\DeclareMathOperator{\Poisson}{Poisson}

\DeclareMathOperator{\T}{\top}
\DeclareMathOperator{\E}{\mathbb{E}}

\newtheorem{theorem}{Theorem}
\newtheorem{corollary}{Corollary}
\newtheorem{lemma}{Lemma}
\newtheorem{definition}{Definition}

\allowdisplaybreaks[3]

\begin{document}
%
\title{Sparse Signal Recovery under Poisson Statistics}

\author{
  \IEEEauthorblockN{Delaram Motamedvaziri, Mohammad H. Rohban, and Venkatesh Saligrama\\}
  \IEEEauthorblockA{
    Electrical and Computer Engineering Department Boston University}}

\maketitle

\begin{abstract}
We are motivated by problems that arise in a number of applications such as Online Marketing and explosives detection, where the observations are usually modeled using Poisson statistics. We model each observation as a Poisson random variable whose mean is a sparse linear superposition of known patterns. Unlike many conventional problems observations here are not identically distributed since they are associated with different sensing modalities. We analyze the performance of a Maximum Likelihood (ML) decoder, which for our Poisson setting involves a non-linear optimization but yet is computationally tractable. We derive fundamental sample complexity bounds for sparse recovery when the measurements are contaminated with Poisson noise. In contrast to the least-squares linear regression setting with Gaussian noise, we observe that in addition to sparsity, the scale of the parameters also fundamentally impacts sample complexity. We introduce a novel notion of Restricted Likelihood Perturbation (RLP), to jointly account for scale and sparsity. We derive sample complexity bounds for $\ell_1$ regularized ML estimators in terms of RLP and further specialize these results for deterministic and random sensing matrix designs.
\end{abstract}

\IEEEpeerreviewmaketitle


\section{Introduction}
\IEEEPARstart{I}{n} this paper, we study the problem of high dimensional sparse model estimation under a Poisson model for observations. This problem is motivated by many practical applications where the observations are the counts of an event. The mean count in these applications depends linearly on a sparse subset of parameters. Our goal is to extract the sparse subset from a potentially large number of parameters. Some of the practical applications motivating our problem include explosive identification based on photon counts in fluoroscopy \cite{cite6}, and eMarketing based on website traffic \cite{cite13}.



We propose a general model that is applicable to a broad class of problems involving Poisson statistics. We consider the case where observations are obtained from heterogeneous sensors or different measurement settings and therefore not identically distributed. To simplify the model, we assume that the rates of the underlying Poisson model for observations are affine functions of some positive signal we want to estimate. In other words, if the signal of interest is $w^*\in \mathbb{R}_+^p$, the $i$-th observation, $y_i$, is distributed as follows:
$$\forall i \in \{ 1,\ldots,n \} : y_i\sim \Poisson(\lambda_{0,i}+a_i^{\T}w^*)$$
where $\lambda_{0,i}$ is the rate of the  background Poisson noise and assumed to be known, and each $a_i=[a_{i,1}, \ldots ,a_{i,p}]^{\T}$ is a distinct vector corresponding to the $i$-th sensor. The collection of these vectors form the sensing matrix, $A=[a_1,\ldots, a_n]^{\T}$. Our goal is to recover the sparse vector, $w^*$, from $\{ y_1, \ldots, y_n \}$.


In this paper we analyze the performance of the $\ell_1$ constrained Maximum Likelihood (ML) decoder. The ML decoder in our Poisson setting is a convex optimization problem involving non-linear objective function. We derive fundamental sample complexity bounds for sparse recovery in the high-dimensional setting. 

In conventional sparse linear least squares regression setting, sample complexity is primarily determined by sparsity for many random sensing matrix designs. While the scale of the parameter vector does influences sample complexity, its impact is somewhat beneficial. Sample complexity improves with scale of the ground truth parameter $w^*$ for a fixed level of noise. In contrast, for our Poisson setting, sample complexity degrades with the scale of the parameter vector. Specifically, sparsity level $k := \lVert w^* \rVert_0$, and parameter amplitude $s := \lVert w^* \rVert_1$ plays a significant role in determining sample complexity. One difference is that the variance of observations grows with $s$ in the Poisson case. A more fundamental reason is that the curvature of the likelihood function decreases with the scale $s$ of the parameter vector. Indeed, for large values of $s$, partial changes $\partial \widehat w = \widehat w - w^*$ in the parameter vector translate to significantly smaller changes of the likelihood function resulting in lack of identifiability. Consequently, unlike the conventional case we inevitably have to suffer the effects of scale for the Poisson case. 


We summarize the main objectives of this paper based on the above discussion. We first characterize sample complexity in terms of constants that account for both sparsity and scale effects of likelihood functions encountered in Poisson type settings. Next by relating these constants to existing eigenvalue based characterizations of sensing matrices we derive sample complexity bounds for different designs. 

We consider the high-dimensional ($p > n$) setting in this paper. High dimensional setting leads to fundamental issues even in the conventional setting. For instance, in least-squares linear regression setting \cite{cite20,Sara} the Hessian is singular. To overcome these issues the loss function is optimized with $\ell_1$ constraints (or the loss function is regularized with an $\ell_1$ penalty) on the parametric space. Following along these lines we also consider optimizing the Poisson likelihood function under $\ell_1$ constraints. These constraints (or penalty) have the effect of constraining the error patterns to a cone of feasible directions (descent cone):
$$
\widehat{w} - w^* \in \mathbb{C} := \{u : \|u_{S^c}\|_1 \leq \|u_S\|_1, |S| \leq k\}
$$
As a result it turns out that we need to ensure that the loss function is ``well-behaved'' only on the feasible cone. A recent important development~\cite{cite20} in this context is to impose strong convexity of the loss function on the feasible cone. While this requirement is generally satisfied for many loss functions including least-squares losses, unfortunately as it turns out, our Poisson case does not satisfy the strong convexity assumption. Strong-convexity of the loss function amounts to the assumption we see a non-trivial change in the loss function as a result of underlying parameter variation regardless of the ground truth $w^*$. This requirement can be viewed as the requirement that the curvature of the loss function is non-vanishing on the cone. In the Poisson case the perturbation in the loss function behaves linearly in large $s$ regimes (i.e. the curvature vanishes in the limit) and so the loss function is no longer strongly convex on the cone.



This issue motivates us to introduce a Restricted Likelihood Perturbation (RLP) constant to characterizes the change in loss function at various amplitudes and sparsity levels in terms of the changes in the parameter vector ($\|\partial \widehat w\| = \epsilon$). We can view RLP implicitly as a condition on sensing matrices. We show that if for some sensing matrix $A$ with bounded elements, the proposed RLP condition is satisfied, then $\ell_1$ constrained ML estimator would converge to $w^*$ with an exponential convergence rate. 

A main drawback of RLP is that it is difficult to directly verify for well-known deterministic and random constructions. A natural condition that has emerged recently in sparse recovery literature is the so called restricted eigenvalue (RE) condition. The RE constant $\gamma_k$ is described as
$$
\frac{1}{n} \| A (\widehat{w} - w^*) \|^2 \geq \gamma_k \| \widehat{w} - w^* \|^2,\,\, \forall \,\,\widehat{w} - w^* \in \mathbb{C}
$$ 
A wide variety of sensing matrices ranging from deterministic to random designs can be shown to satisfy this requirement. 

This motivates us to express RLP constant in terms of the RE parameter $\gamma_k$. We then obtain an expression that relates RLP constant as a function of $\gamma_k$ and amplitude $s$. This leads to an explicit expression for sample complexity in terms of scale and sparsity making the dependence transparent. It follows that for a fixed scale the sample complexity is completely determined by sparsity and is simply a function of the sensing matrix. This dependence is generally similar to that obtained for CS for bounded positive sensing matrices. On the other hand for a fixed sparsity level our results shows that sample complexity has an inverse relationship to scale.

In particular our results apply to both deterministic and random sensing matrices and we present several results for both cases. We also conduct several synthetic and real-world experiments and demonstrate the  tightness of the oracle bounds on error as well as the efficacy of our method. Specifically, it has been suggested in the literature that LASSO can handle exponential family noise such as that arises in our application~\cite{cite5}. It turns out that $\ell_1$ constrained ML estimator uniformly outperforms LASSO and has significantly superior performance in many interesting regimes.

The paper is organized as follows: In Section \ref{PS}, we introduce the notation and state our sparse estimation problem. Section \ref{OA} describes our theoretical results on the convergence of the regularized ML decoder. The proofs are briefly sketched in Section \ref{PrfSk}. The numerical results for different interesting scenarios are demonstrated in Section \ref{NR}. Finally, the detailed proof of the main theorems and lemmas are provided in Section \ref{AP}.

\subsection{Related Work}

Parameter estimation for non-identical Poisson distributions has been studied in the context of Generalized Linear Models (GLMs). However, our model is inherently different from the exponential family of GLM models that has been studied in \cite{cite10,cite12,cite14,cite15}. In particular the GLM model corresponding to the Poisson distributed data studied in the literature has the following form:
\begin{align*}
\text{Model I : } & ~~~~ \Pr(y_i=k)=\Poisson(\exp\left(a_i^{\T}w\right))\\
&\propto \exp\left(k\left(a_i^{\T}w\right)\right)\exp\left(-\exp\left(a_i^{\T}w\right)\right)
\end{align*}
Therefore, the log likelihood has the form:
\begin{align*}
\mathcal{L}_1(w)=\sum_{i=1}^ny_i\left(a_i^{\T}w\right)-\exp\left(a_i^{\T}w\right)
\end{align*}
In contrast, in the setting we are interested in, the observations are modeled as follows:  
\begin{align*}
\text{Model II : }\Pr&(y_i=k)=\Poisson(\lambda_{0,i}+a_i^{\T}w )\\
&\propto \left(\lambda_{0,i}+a_i^{\T}w\right)^k \exp\left(-\left(\lambda_{0,i}+a_i^{\T}w\right)\right)
\end{align*}
and the log likelihood function has the form:
\begin{align*}
\mathcal{L}_2(w)=\sum_{i=1}^ny_i\log \left(\lambda_{0,i}+a_i^{\T}w\right)-\left(\lambda_{0,i}+a_i^{\T}w\right)
\end{align*}

As a statistical model there are several differences between the two models. We observe that imposing sparsity on $w$ in Model I corresponds to smaller number of multiplicative terms. On the other hand, $w$ being sparse in Model II results in fewer number of additive terms in the Poisson rate of the corresponding model. 
At a more fundamental level the loss function (negative log-likelihood) for Model I has an exponential term ($\exp\left(a_i^{\T}w\right)$). The assumptions of strong convexity on the feasible cone are readily satisfied. Consequently, unlike our case, the issue of signal amplitude no longer arises for this model. 
%
%
%
Therefore, we can view this model as an instance of a general class of sparse problems. Indeed, \cite{cite14} studies the convergence behavior of $\ell_1$ regularized ML for exponential family distributions and GLM in this context. The bounds on error for sparse recovery of the parameter are based on the RE condition. Moreover, in order to get useful bounds on estimation error of GLM, they additionally need the natural sufficient statistic of the exponential family to be sub-gaussian. This condition could clearly be violated in our setting where the data is Poisson distributed and there is no constraint on the sensing matrix to be sub-gaussian.

More generally  \cite{cite10} describes a unified framework for  analysis of regularized $M-$ estimators in high dimensions. They also mention extension of their framework to GLMs and describe ``strong convexity" of the objective function as a sufficient condition to obtain consistency of M-estimators under Model I. As we described this requirement of strong-convexity is not consistent with our model. In addition the statistical aspects in that work requires that the components of the sensing matrix be characterized by \textit{sub-Gaussian distributions}, which we do not require here. 



Statistical guarantees for sparse recovery in settings similar to model II have been provided in \cite{cite11,cite112,cite113} in the context of photon limited measurements. They assume that the observations are distributed as follows  
$$y_i\sim \Poisson(a_i^{\T}w^*)$$
where elements of the signal $w^*$ and sensing matrix are positive, and the sensing matrix satisfies the so-called Flux Preserving assumption:
$$ \sum_{i=1}^{n} (A w^*)_i \leq \sum_{i = 1}^{n} w^*_i .$$

The latter assumption arises in some photon counting applications, like imaging under Poisson noise, where the total number of expected measured photons cannot be larger than the intensity of the original signal. The upper bound on reconstruction error of the constrained ML estimator is given in the paper \cite{cite112}. Surprisingly, the upper bound scales linearly with the number of measurements. However, this sounds reasonable under the Flux Preserving assumption. In fact this behavior is due to the fact that for a fixed signal intensity, more measurements lead to lower SNR for each observation. As a result, unlike conventional compressive sensing bounds, the estimates do not converge to the ground truth with increasing the sample size.
Nevertheless, Flux Preserving constraint does not arise in our setting and consequently the application and methods of analysis are different.

In summary the fundamental differences in the underlying model as well as the different assumptions in the sensing matrices from the previous work warrants new analysis techniques which is the subject of this paper.  



\subsection{Applications}

In the sequel, we will introduce two applications, which motivate the model described earlier: \\ 

{\bf 1. Explosive Identification:} In the explosive identification example, an unknown mixture of explosives is exposed to different  fluorophores. The goal is then to estimate the mixture components based on the observed fluorophores photon counts. Here $y_i$ and $\lambda_{0,i}$ could be considered as the photon counts and background emission rates for fluorophore $i$. $\lambda_{0,i}$ is measured before the exposure of the explosive mixture and is known. $a_{ij} < 0$ represents a quenching effect of explosive $j$ on fluorophore $i$, and $w^*_j$ is the weight of explosive $j$ in the mixture~\cite{cite6}. \\

\begin{figure}
\begin{minipage}[b]{1\linewidth}
  \centering
  \centerline{\includegraphics[width= 8cm]{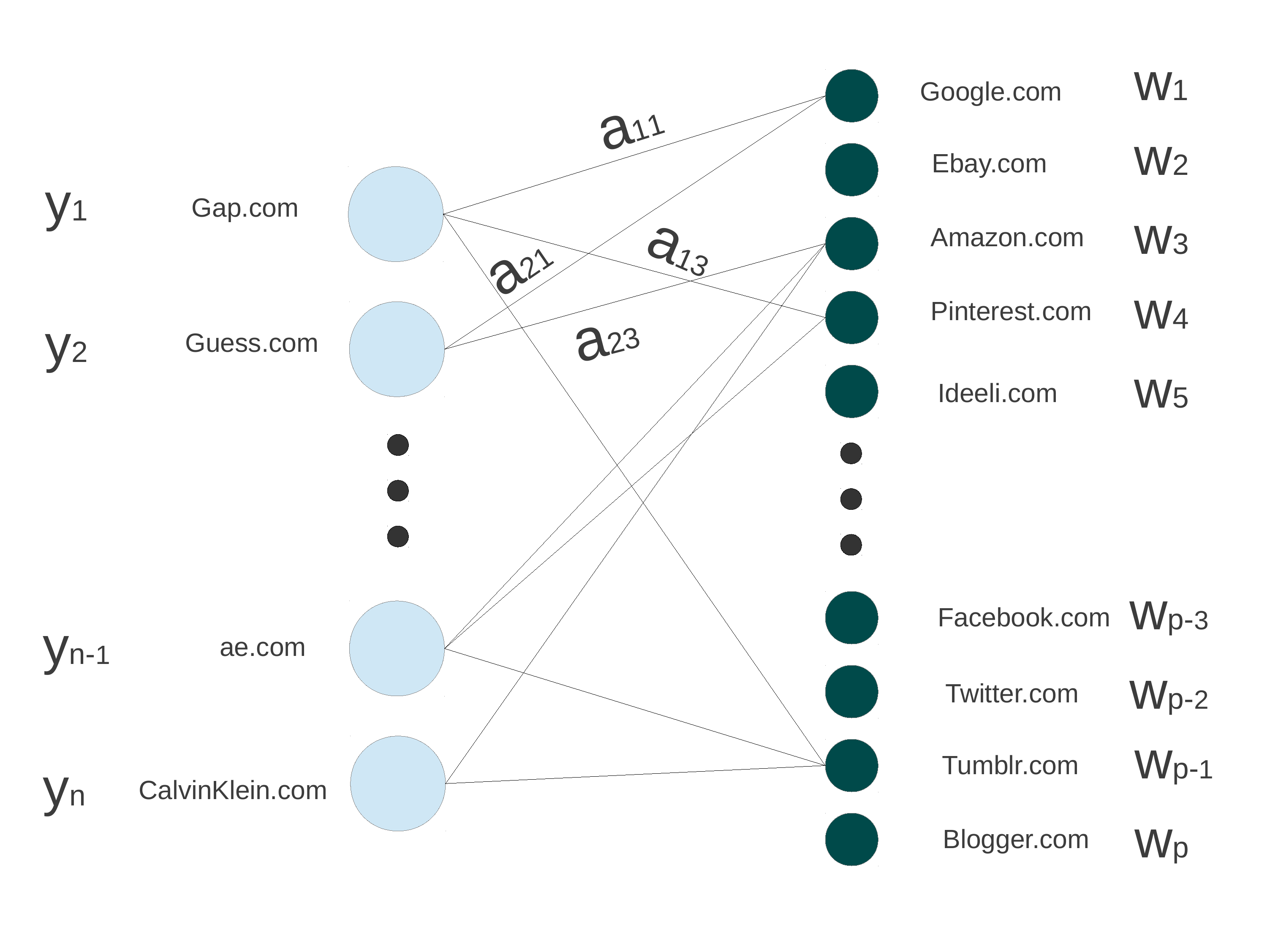}}
\caption{ Our eMarketing model: right nodes are the business websites, left nodes are the advertisement websites. A connecting edge, $a_{ij}$, is the number of backward links purchased by the business website $i$ from the advertisement website $j$.
}
\label{MarketingFig}
\end{minipage}
\end{figure}

{\bf 2. eMarketing:} In the eMarketing example~\cite{cite13}, weekly traffic of different websites within the same market are measured. The traffic of site $i$, denoted by $y_i$, is assumed to be affected by the number of bought links to it on different advertisement websites. Each advertisement website $j$ is also assumed to contribute to the visiting rate by a fixed weight $w^*_j$. These weights could be viewed as measures of popularity/dominance of advertisement website $j$. Moreover, $\lambda_{0,i}$ is the average traffic that visits website $i$ directly (not through intermediate advertisement website) and is acquired through online statistics of the website. $a_{i,j} > 0$ is the number of backward links that business website $i$ has bought from advertisement website $j$. Estimation of $w^*_j$'s can lead to discovery of dominant advertisement websites in some industry. Fig. \ref{MarketingFig} illustrates this eMarketing model.

\section{Problem Setup}
\label{PS}
\subsection{Notation}
\label{not}
We will use the following notation in this paper. First, we assume that $w^* \in \mathbb{R}^{p}_{+}$ is the vector of true parameters that generates Model II. For ease of exposition we assume elements of $A$ to be positive, although similar bounds could be obtained for the case that $a_{i,j}$'s are negative. Let $a_i$'s be the rows of the sensing matrix $A$, and consider the following simplifying notations for later use: 

\begin{itemize}
\item $a_{max} :=\max_{i,j} a_{i,j}$
\item $\lambda_{w^*,i} := \lambda_{0,i}+a_i^{\T}w^*$
\item $\lambda_{0} := \max_i\lambda_{0,i}$
\item $k := \lVert w^* \rVert_0,  ~~~ s := \lVert w^* \rVert_1$
\item $\lambda_{max} := \max_{i, w^*} \lambda_{w^*,i} = \lambda_0 + sa_{max}$
\item $\lambda_{min} := \min_i\lambda_{w^*,i}$
\end{itemize}
We consider $n$ independent Poisson distributed observations generated as:
$$\forall i \in \{ 1,\ldots,n \} : y_i \sim \Poisson(\lambda_{w^*,i})$$
This model arises in applications where the measurements are superposition of independent arrival processes of interest contaminated by some independent background arrival.

We define $Q$ as the normalized negative log-likelihood for these observations :
\begin{equation} \label{Qfunc}
Q(w) :=  -\frac{1}{n} \sum_{i=1}^{n} y_i \log(\lambda_{0, i} + a_i^{\top}w) - a_i^{\top} w 
\end{equation}
Moreover we define $\overline{Q}(w)$ as the expected value of $Q(w)$ :
\begin{equation} \label{Qbarfunc}
\overline{Q}(w) := \mathbb{E}(Q(w)) =  -\frac{1}{n} \sum_{i=1}^{n} (\lambda_{w^*,i}) \log(\lambda_{0, i} + a_i^{\top}w) - a_i^{\top} w 
\end{equation}

Finally, note that all probabilities are calculated conditional on the true parameters vector, $w^*$. In other words, for any event $A$:
\begin{equation} \label{condD}
\Pr\{A\} := \Pr\{A|w^*\}
\end{equation}
\subsection{Problem Formulation}

Our goal is to recover $k$-sparse weight vector, $w^*$, from the observations $y_i$'s, the sensing matrix $A$, and the background noise rates $\lambda_{0,i}$'s. Although $y_i$'s are non-identically distributed, their distributions are related through $w^*$. Hence, estimating the weight vector $w^*$ can be interpreted as a parameter estimation problem using $n$ independent non-identical Poisson distributed samples, which are related through $k$ non-zero elements of $w^*$. 

We study the high dimensional problem where the number of parameters $p$ can grow rapidly with $n$, and $k$ can scale with $p$. Our goal is to prove that under appropriate conditions on $a_i$'s, $\widehat{w}$, the $\ell_1$ regularized ML estimate of $w^*$ from $y_i$'s, is consistent with the ground truth:
$$\lim\limits_{n \rightarrow \infty}\Pr\{\lVert\widehat{w}-w^*\rVert_2\geq\epsilon\}=0$$
where this probability is conditional on $w^*$ as described in Eqn. $(\ref{condD})$. Moreover, we want to show exponential rate of convergence with respect to the number of observations:
$$\Pr\{\lVert\widehat{w}-w^*\rVert_2\geq\epsilon\}\leq C'\exp({-nC})$$
where $C$ and $C'$ are some positive constants.

\subsection{Regularized Maximum Likelihood}
We first describe sparse sets and amplitude-constrained sets for later use.
Let  $\Gamma_k$ be the set of k-sparse signals:  $$\Gamma_k =\{w|\,w\geq 0, \lVert w\rVert_0\leq k \}$$ 
Let $\Theta_s$ denote amplitude-constrained sets of scale $s$: $$\Theta_s =\{w|\,w\geq 0, \sum_{j=1}^p w\leq s \}$$ 


Note that in our problem, the observations are independent and follow a Poisson distribution. With this in mind we propose the following constrained ML estimation:
\begin{equation}
\label{eq2}
\widehat{w}=\arg\max_{{w\in \Theta_s}} Q(w) \stackrel{\Delta}{=} \arg\max\limits_{\substack{w\in \Theta_s}}\,\log p(y_1, \ldots ,y_n|{w})
\end{equation}
The constrained maximization problem defined in Eqn. (\ref{eq2}) is equivalent to the following  minimization problem for a suitably chosen $\eta_s$, through a Lagrangian formulation:
\begin{equation}
\label{eq3}
\widehat{w} = \argmin_{i : w_i \geq 0} Q(w) +\eta_s\sum_{j=1}^p w_j
\end{equation}
These problems are convex and can be solved efficiently by conventional optimization algorithms to find the global optimum. It needs to be mentioned that $y_i$'s are not identically distributed so the consistency of the resulting estimates does not trivially follow from the consistency of ordinary maximum likelihood. In the next section, we will describe sufficient conditions on consistency of regularized ML estimation.

\section{Main Results}
\label{OA}
We will first abstractly state general results for deterministic sensing matrices. We will then relate these results to specific sensing designs. In this context we introduce a new notion that jointly accounts for scale as well as sparsity.  
\subsection{Restricted Likelihood Perturbation}
We introduce the notion of Restricted Likelihood Perturbation to jointly account for scale and sparsity. 
\begin{definition} {Restricted Likelihood Perturbation, $RLP(A, \beta_{s,k})$ :}
Let $\delta_{s,k}(\epsilon)$ denote:
\begin{equation} \label{deltask}
\delta_{s, k}(\epsilon) := \min\limits_{\stackrel{\lVert w - w^* \rVert_2 = \epsilon}{w \in \Theta_s, ~ w^* \in \Gamma_k \cap \Theta_s}} \overline{Q}(w) - \overline{Q}(w^*)
\end{equation}
The sensing matrix $A$ is said to satisfy RLP condition, when for each $s$ and $k$, there exist a constant $\beta_{s,k}$ such that \hbox{$\delta_{s,k}(\epsilon)/\epsilon^2 \geq \beta_{s,k}$.}
\end{definition}
Note that $\beta_{s,k}$ is the minimum perturbation rate in the averaged loss function $\overline{Q}$ caused by the change in the parameter value ${\widehat{w} - w^*}$.
The following result characterizes the estimation error in terms of $\beta_{s,k}$. 
\begin{theorem} \label{mainTh}
Suppose that $A$ satisfies $RLP(A, \beta_{s,k})$, 
\hbox{$\lambda_{w^*, i} \geq \lambda_{min}$,}  $\lambda_0 = \mathcal{O}(1)$, $a_{i,j} \leq a_{max} = \mathcal{O}(1)$ for all $i$ and $j$, and $e$ be a real number with $0 < e < 1$.
Suppose further that the number of measurements $n$ satisfies
$$ n \geq  \frac{c_1 s \log^2(s) \log\frac{2}{e} }{\beta_{s,k}^2 \epsilon^4} $$
and $\epsilon$ is small enough such that
$$ 0 < \epsilon \leq \sqrt{\frac{c_2 \lambda_{min} \log s}{\beta_{s,k} \max(c^{\prime}, \sqrt{s})}}. $$
Then the probability of error for the constrained ML estimate of Eq.~\ref{eq2} being greater than $\epsilon$ can be bounded as follows  :
$$ \Pr\{ \lVert \widehat{w} - w^* \rVert_2 \geq \epsilon \} \leq e$$
where
$c_1$, $c_2$, and $c^{\prime}$ are universal constants, which do not depend on $A$ or ${w}$.
\end{theorem}
We note that if the background rate $\lambda_{\min}$ approaches zero the admissible $\epsilon$'s approach zero as well. In the limit this implies that the sample complexity approaches infinity. This makes sense because as $\lambda_{\min}$ approaches zero, we lose identifiability in the limit. \\

{\bf Remark}:
Theorem \ref{mainTh} suggests that the sample complexity hinges mainly on two factors: $s$ and $\beta_{s, k}$. The effect of amplitude $s$ could be explained intuitively as follows. For $\widehat{w}$ to be within a close distance of $w^*$ with high probability, $Q$ should converge in some notion to $\overline{Q}$. Since the variances of the Poisson distributed observations are equal to the rates $\lambda_{w^*, i}$, more samples are required for $Q$ to concentrate around $\overline{Q}$ for large $\lambda_{max} = \mathcal{O}(s)$. This effect does not appear in conventional compressive sensing where Gaussian distributed samples are used. 

Note that $\beta_{s,k}$ conceptually measures the curvature of the objective function $\overline{Q}$ around $w^*$. A small curvature is detrimental in that $Q$ would be more tightly concentrated around $\overline{Q}$ for a given $\|\partial \widehat{w}\| = \|\widehat{w} -w^*\| =\epsilon$. Consequently small curvature degrades sample complexity.

$\beta_{s,k}$ could also be viewed as a measure of {\it identifiability} of $w^*$. If ${\beta_{s,k}}$ is zero or close to zero, there would be no unique solution to the likelihood maximization problem, and hence there would be no hope for $\widehat{w}$ to converge to $w^*$ in the worst case. $\square$\\

\noindent
{\it Relationship to Restricted Eigenvalue Condition:} \\
Although RLP precisely accounts for sample complexity, it is hard to verify for well-known deterministic and random constructions. On the other hand, there has been a significant amount of literature on the so called Restricted Eigenvalue condition (RE) for sensing matrices. RE is the basis for analysis of many noisy sparse recovery methods. Our goal is to find a relationship between RLP and RE here so that our sample complexity results are transparent for many of these constructions.
%
\begin{definition} {Restricted Eigenvalue condition, RE$(A, \gamma_k)$:} 
There exists a constant $\gamma_k > 0$, such that for any set of indices $S$ satisfying $|S| \leq k$, and vector 
$$u \in \mathbb{C}(S) := \{ u \neq {\mathbf 0} : \lVert u_S \rVert_1 \geq \lVert u_{S^c} \rVert_1 \},$$
we have:
\begin{equation} \label{REC}
 \frac{1}{n} \|Au\|_2^2\geq \gamma_k \|u\|^2_2
\end{equation}
where $u_S$ is the restriction of the vector $u$ to the indices in $S$, and $S^c = \{1, \ldots, p \} \setminus S$.
\end{definition}

{\bf Remark} : The $1/n$ factor on the left hand side of Eqn. \eqref{REC} can be considered as a column normalization of $A$, i.e. each column is divided by $\sqrt{n}$.

RE condition is a well known sufficient condition for consistency of several sparse recovery algorithms. Specifically, its various forms were used to derive the oracle inequalities for LASSO and Dantzig selector \cite{cite20}, \cite{Van09}. 

There are a number of well known results for random designs $A$, for which Definition 1 holds with high probability in terms of $n$ \cite{Rask10}. For example, consider the case that elements of $A$ are i.i.d. samples from a subgaussian distribution. Then, Definition 1 is satisfied for all $n \geq c k\log(p)$, with probability at least $1 - c_1 \exp(-c_2 n)$, where $c$, $c_1$, and $c_2$ are universal constants \cite{ShuZh09}. Moreover, in these cases $\gamma_k$ is invariant to sparsity level $k$, so long as $n \geq c k\log(p)$. 

In our experimental data, however, matrix $A$ is deterministic and given to us. In general, testing RE condition is an NP-hard problem. Nevertheless, our numerical results still show fast rate of convergence for regularized ML. $\square$



The following result characterizes the relationship between RLP and RE conditions.
\begin{theorem}
\label{mainTh2}
Suppose that RE$(A, \gamma_k)$ holds, elements of $A$ are bounded, and $\lambda_0 = \mathcal{O}(1)$. Then, $A$ satisfies $RLP(A, \beta_{s,k})$ for the following value of $\beta_{s, k}$:
\begin{equation*}
\label{Delt} 
\beta_{s,k} = \frac{c\gamma_k}{s}
\end{equation*}
where $c$ is a universal constant.
\end{theorem}
{\bf Remark:} Note that the impact of the sparsity and scale parameter $s$ is now transparent. Indeed sparsity and scale are separately captured by the expression. It follows that for a fixed scale the sample complexity is completely determined by sparsity and is simply a function of the sensing matrix. On the other hand for a fixed sparsity level the constant $\beta_{s,k}$ decreases with scale $s$.

\begin{corollary} \label{Cor1}
Under the conditions of Theorems \ref{mainTh} and \ref{mainTh2}, if 
$$ n \geq  \frac{c_1 s^3 \log^2(s) \log\frac{2}{e} }{\gamma_k^2 \epsilon^4} $$
and,
$$ 0 < \epsilon \leq \sqrt{\frac{c^{\prime} \lambda_{min} s 
\log s}{\max(c_2, \sqrt{s}) \gamma_k}} $$
then
$$ \Pr\{ \lVert \widehat{w} - w^* \rVert_2 \geq \epsilon \} \leq e$$
where $c_1$, $c_2$, and $c^{\prime}$ are universal constants and $0<e<1$.
\end{corollary}
{\bf Remark:}
Theorem \ref{mainTh2} suggests that under the RE condition, $\beta_{s,k}$ decreases with the increase in the signal amplitude $\lambda_{max} = \mathcal{O}(s)$.  This can be attributed to the fact that $\overline{Q}(w) - \overline{Q}(w^*)$ as function of $X_i := (w - w^*)^{\top} a_i$ behaves sub-quadratically (and in fact almost linearly) for large $X_i$. In fact, the perturbation in $\overline{Q}$ could be shown to scale  in proportion to $-\lambda_{w^*, i} \log(1+X_i/\lambda_{w^*, i})+X_i$. However, RE condition lower bounds a quadratic function of $X_i$'s. Therefore, if RE condition is to be used to lower bound the perturbation rate in $\overline{Q}$, an extra factor involving the signal $\ell_1$ intensity would inevitably appear in the lower bound of perturbation. This factor does not arise in conventional compressive sensing as the objective functions there are quadratic in $X_i$'s. 

\begin{figure}
\begin{minipage}[b]{1\linewidth}
  \centering
  \centerline{\includegraphics[width= 8cm]{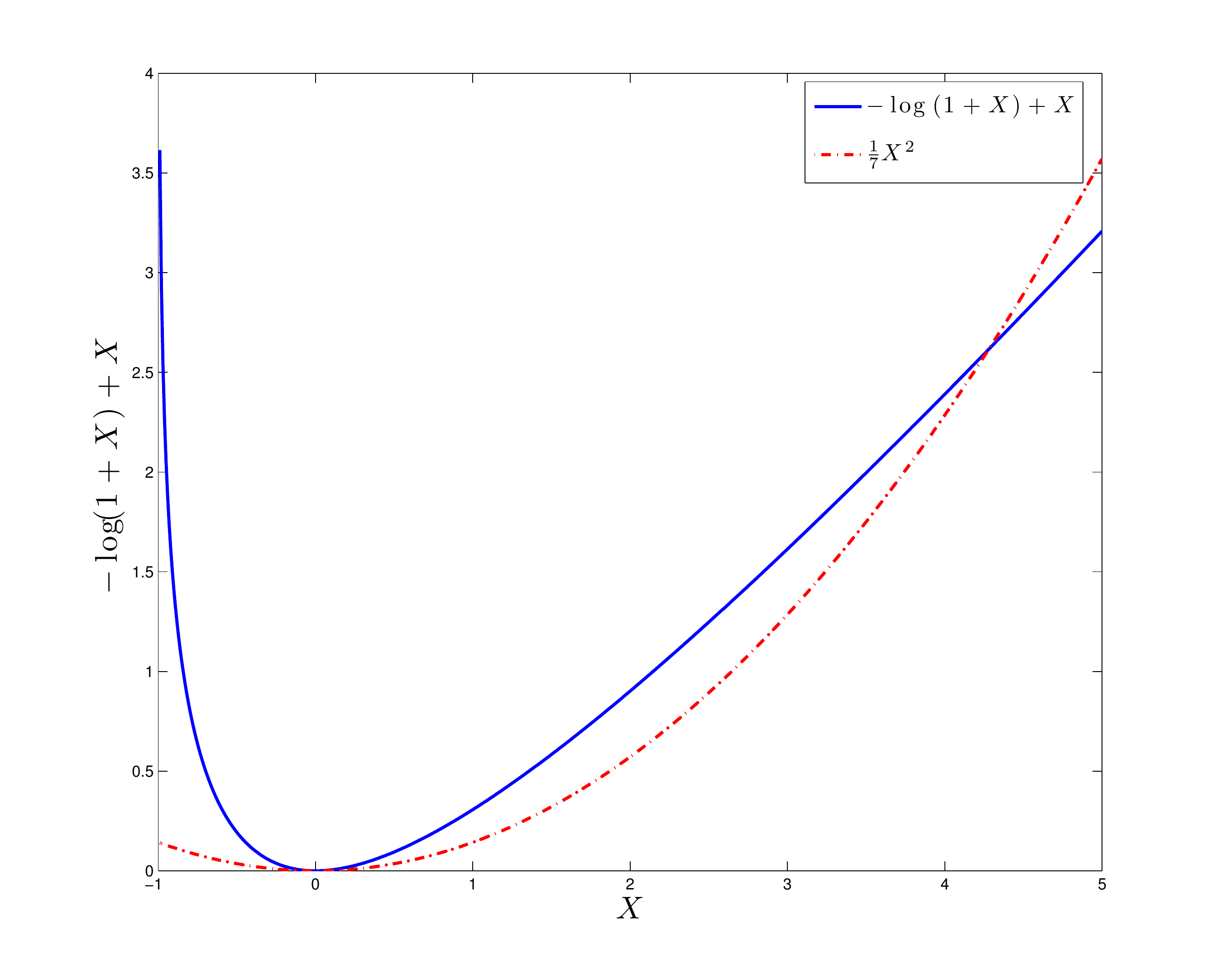}}
\caption{$\overline{Q}$ behaves like $-\log(1+X)+X$. It scales almost linearly for large X. Therefore, to be lower bounded by a quadratic function, an extra scale factor is needed. This factor would be inverse of the maximum value that $X$ can take, which would be proportional to $1/s$ in our problem setting.}
\label{LossComp}
\end{minipage}
\end{figure}

\begin{figure}
\begin{minipage}[b]{1\linewidth}
  \centering
  \centerline{\includegraphics[width= 8cm]{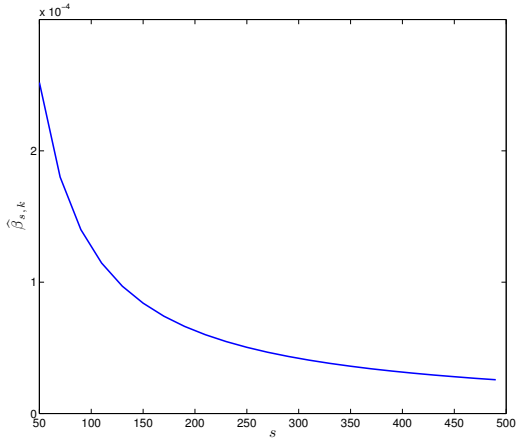}}
\caption{$\widehat{\beta}_{s,k}$ is an approximation to $\beta_{s,k}$ for $p = 200$, $k = 10$, and $n = 40$. Entries of $A$ are i.i.d. instantiations of the uniform distribution between 0 and 1. The plot shows that $\widehat{\beta}_{s,k}$ decreases with the increase of parameter amplitude $s$.}
\label{ApproxBet}
\end{minipage}
\end{figure}

We illustrate this point in Fig. \ref{LossComp}. The plot depicts perturbation in $\overline{Q}$ and compares it to a quadratic loss function. We see that $\overline{Q}$ is essentially linear for large $s$. This issue can also be seen from Fig. \ref{ApproxBet}. Here an approximation of $\beta_{s,k}$ is obtained by taking the minimum of perturbation rates in $\overline{Q}$ over $10^6$ random points in the cone of feasible directions. Fig. \ref{ApproxBet} clearly demonstrates the inverse relationship between $\beta_{s,k}$ and $s$. $\square$

The result in the last corollary gives the number of required samples for a certain level of accuracy. Alternatively, it is possible to restate the result to find a high probability oracle on the error:
\begin{corollary} \label{Cor2}
Under the conditions of Theorems \ref{mainTh} and \ref{mainTh2}, if 
\begin{equation} \label{NBndNew}
n \geq \left(\frac{C_0 s}{\lambda_{min}} \right)^{\frac{5}{2}}
\end{equation}
the $\ell_2$ error of regularized maximum likelihood estimation follows: 
\begin{equation} \label{ErrBndNew}
\lVert \widehat{w} - w^* \rVert_2 \leq \frac{C^{\prime} s^{\frac{3}{4}} \log^{\frac{1}{2}}(s) }{\gamma_k^{\frac{1}{2}} n^{\frac{1}{5}}} 
\end{equation}
with probability of at least $1-C_1\exp(- C_2 n^{\frac{1}{5}})$, where $C_0$, $C_1$, $C_2$, and $C^{\prime}$ are constants. 
\end{corollary}

\noindent
{\it Rate of Convergence:} Unlike traditional compressed sensing where error scales with $n^{-\frac{1}{2}}$, the error here scales with $n^{-\frac{1}{5}}$. We have experimentally verified the tightness of this result in section \ref{TRV}. \\
{\it Impact of Background $\lambda_{0,i}$:} For the sake of exposition we have assumed $\lambda_{0, i}$'s are bounded. When $\lambda_{0,i}$'s are variable we would have to replace $s$ by $\lambda_{max} :=\lambda_0 + a_{max} s$ in all the equations, with $\lambda_0 := \max_{i} \lambda_{0, i}$ and $a_{max} := \max_{i,j} A_{i,j}$. 
Assuming that $\lambda_{0, i}$'s are equal, the minimum Poisson rate changes with $\lambda_0$, i.e. $\lambda_{min} \geq \lambda_0$. The impact of $\lambda_0$ on the error is subtle. On the one hand for large $\lambda_0$, the lower bound on $n$ in Eq. \eqref{NBndNew} decreases leading to a wider ranger of admissible sample sizes. On the other hand, the upper bound on error in Eq. \eqref{ErrBndNew} scales as $\lambda_{max}^{3/4} \log^{1/2} \lambda_{max}$, larger $\lambda_0$ results in a larger error bound. Still when the sample size is sufficiently small, the first effect dominates the second one, i.e. large $\lambda_0$ appears to be beneficial. For large sample size estimation error generally increases with $\lambda_0$.

\subsection{Random Design}

It has been so far assumed that the measurement matrix $A$ satisfies the RE condition. The next theorem gives the error bound for random matrix constructions with bounded and positive elements. 

\begin{theorem} \label{mainTh4}
If elements of $A$ are i.i.d. samples from a distribution with bounded support on $\mathbb{R}_{+}$, then for a suitable choice of $\Theta_s$ (or equivalently, if $\eta_s$ is chosen appropriately), any 
\begin{multline*}
n\geq \max\Bigg( c_1 k^2 \log ( p ) \log^3 (c_2 k \log( p )) \log\left(\frac{1}{e}\right),   \\\frac{c^{\prime} k^2 s^3 \log^2(s) \log\frac{4}{e} }{\epsilon^4} \Bigg)
\end{multline*}
and
$$ 0 < \epsilon \leq \sqrt{\frac{c_4 k^2\lambda_{min} s 
\log s}{\max(c_3, \sqrt{s})}},$$
we have:
\begin{equation*}
\Pr\{\|\widehat{w}-w^*\|_2\geq\epsilon\} \leq e
\end{equation*}
where $c^{\prime}$, $c_1, \ldots, c_4$ are constants. 
\end{theorem}

\noindent
{\bf Remark:}
We observe that unlike the traditional settings the sample complexity grows primarily with $k^2$ for this particular sensing matrix. The main reason that underlies this effect is the requirement that $A$'s elements have to be supported on a set of positive numbers. This constraint precludes the possibility of sub-gaussian constructions of $A$. However, most of the known random design constructions with linear measurement complexity $\mathcal{O}(k)$ require $A$'s element to be at least sub-gaussian for RE condition to be satisfied. If $A$'s elements are not sub-gaussian, the best known result requires at least $\mathcal{O}(k^2)$ measurements to satisfy RE condition with high probability \cite{Rud10}. Therefore, if the positivity requirement is eliminated, the usual $\mathcal{O}(k)$ measurement complexity would be possible.

\section{Proof Sketch} \label{PrfSk}
\subsection{Restricted Likelihood Perturbation}
We use the idea of {\it Extremum Estimators} in our proof. These are a broad class of estimators for parametric models calculated through maximization (or minimization) of an objective function $Q(w)$, which depends on the data \cite{cite17}.  

\begin{lemma} \label{upperBnd}
If $\widehat{w}$ and $w^*$ are the minimizers of $Q(w)$ and $\overline{Q}(w)$ respectively subject to $w \in \Theta_s$, it follows that:
\begin{align} \label{unifConv}
\Pr\{\lVert w^*- \widehat{w}\rVert_2 \geq \epsilon \}
 \leq  \Pr\left\{ \sup_{w \in \Theta_s}  | Q(w)-\overline{Q}(w) |  \geq \frac{\delta_{s,k}(\epsilon)}{2} \right\}
\end{align}
where $\delta_{s,k}(\epsilon)$ is defined in Eq. \eqref{deltask}
\end{lemma}
\begin{proof}
The detailed proof is provided in Appendix A.
\end{proof}
Intuitively, $\delta_{s,k}$ in Lemma \ref{upperBnd} represents the minimum increase in the function $\overline{Q}(\widehat{w})$, when $\widehat{w} \in \Theta_s$ is $\epsilon$ far away from the function minimizer $w^{*}$. When the function $\overline{Q}(\widehat{w}) - \overline{Q}(w^*)$ is strongly convex in terms of $u := \widehat{w}-w^*$, $\delta_{s,k}$ would be strictly positive, for $\epsilon > 0$. Strong convexity is a strong condition in high dimensional settings. However, in our setting $u$ could be shown to belong to a feasible cone (see Section~1). Therefore we require strong convexity only on a {\it restricted} set of directions for a fixed scale $s$. Assuming that $\delta_{s,k} > 0$, if $Q$ is uniformly convergent to $\overline{Q}$ on $\Theta_s$, the right hand side of \eqref{unifConv} will converge to zero for a fixed $s$. Then, this will imply the consistency of $\widehat{w}$.

It can be seen that $w^{*}$ is the minimizer of $\overline{Q}$ over $\Theta_s$. Hence if $\delta_{s,k} > 0$, preconditions of Lemma \ref{upperBnd} are satisfied and we may use the upper bound in Eqn. \eqref{unifConv}.

Next, we are going to upper bound the right hand side of Eqn. \eqref{unifConv} using the following lemma :
\begin{lemma} \label{UltLem}
For $0 < \delta \leq \frac{c_1 \lambda_{min} \log(\lambda_{max})}{\max(c^{\prime}, \sqrt{s})}$ and $\lambda_0 = \mathcal{O}(1)$, $Q$ is uniformly concentrated around $\overline{Q}$:
\begin{equation}
\label{UltEq}
\Pr\left\{\sup_{w\in\Theta_s} |Q(w) - \overline{Q}(w)| \geq \frac{\delta}{2}\right\} \leq \exp{\left(\frac{-c_2 n\delta^2}{s \log^2 s} \right)}
\end{equation}
where $c_1$, $c_2$, and $c^{\prime}$ are two universal constants.
\end{lemma}
\begin{proof}
The detailed proof is provided in Appendix A.
\end{proof} 
Theorem \ref{mainTh} could be proved by first noting that according to $RLP(A, \beta_{s,k})$, $\delta_{s,k}(\epsilon)$ could be replaced in all of the bounds in two last lemmas by $\beta_{s,k} \epsilon^2$, and then combining the results of the two lemmas.
To prove Theorem \ref{mainTh2}, we approximate the perturbation in $\overline{Q}$ by its truncated taylor series in terms of $a_i^{\top} u$, where $u = \widehat{w} - w$. As elements of $A$ are assumed to be bounded, $a_i^{\top} u$ would be bounded too, and the perturbation would be lower bounded by a quadratic function in $a_i^{\top} u$'s. Then, we use RE condition to lower bound the quadratic function in the last step. The detailed proof is given in the Appendix A.

Corollary \ref{Cor1} is the direct consequence of combining Theorems \ref{mainTh} and \ref{mainTh2}. Corollary \ref{Cor2} is obtained by first rewriting Corollary \ref{Cor1} as:
$$ \Pr(\lVert \widehat{w} - w \rVert_2 \geq \epsilon) \leq 2 \exp\left( \frac{-c n \epsilon^4 \gamma_k^2}{s^3 \log^2 s} \right) $$
Then, by setting $\epsilon = {C^{\prime} s^{\frac{3}{4}} \log^{\frac{1}{2}} }/({\gamma_k^{\frac{1}{2}} n^{\frac{1}{5}}})$, the desired bound on the $\ell_2$ error would be obtained by a probability of at least $1 - C_1 \exp(-C_2 n^{\frac{1}{5}})$.



\subsection{Random Design}
As Theorem \ref{mainTh2} requires RE condition as well as elements of $A$ to be bounded, one may extend these results to case that $A$ is chosen 
randomly and the preconditions are satisfied with high probability. RE condition can be guaranteed with high probability for various classes of random designs. Specifically, if rows of $A$
are drawn from an ensemble of isotropic subgaussian random variables (or a linear transformation of them), it is well known that RE condition will be satisfied with overwhelming probability ($1 - 2 \exp(-c n) $) when $n = \Omega(k \log( p ))$ with $\gamma_k = \Omega(1)$\cite{ShuZh09}. 

However, our setting needs all entries of $A$ to be positive, which is not satisfied for a sub-gaussian ensemble. Therefore, we use the following Lemma from \cite{Rud10}, which guarantees RE condition for the case that elements of $A$ are i.i.d. samples from a bounded random variable (a variation of Theorem 1.8 to adapt to our definition of RE condition and notations).
\begin{lemma} \label{RERa}
If elements of $A$ are i.i.d. samples from a distribution with  bounded support, and \hbox{$n = \Omega(k^2 \log ( p ) \log^3 (k \log ( p )))$}, then $A$ satisfies RE condition with $\gamma_k = \Omega(\frac{1}{k})$ with probability at least $1 - \exp(-c n/k)$, where $c$ is a constant.
\end{lemma}

Based on Lemma \ref{RERa}, Theorem \ref{mainTh4} can be obtained from Corollary \ref{Cor1}, by conditioning the probability of error on the fact that $A$ satisfies RE with $\gamma_k = \Omega(1/k)$. Then, applying a union bound on the probability of the previous event and the event that $A$ satisfies RE, the convergence rate of error in Theorem \ref{mainTh4} would be obtained.


\section{Numerical Results}
\label{NR}
\subsection{Tightness of the Error Bounds} 
\label{TRV}
\begin{figure}
\begin{minipage}[b]{1\linewidth}
  \centering
  \centerline{\includegraphics[width= 9cm]{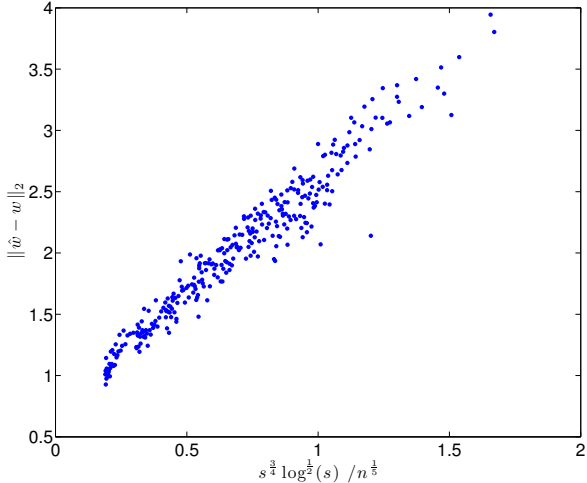}}
\caption{$\ell_2$ error against the upper bound for the error in Corollary \ref{Cor2} for parameter amplitude values $s$ and sample sizes $n$. Evidently, our derived bound of Eq.~\ref{ErrBndNew} is tight as the points show almost a linear relationship between error and the upper bound.}
\label{Cor2Fig}
\end{minipage}
\end{figure}

In this section, our goal is to experimentally verify the importance of amplitude effect described in Corollary \ref{Cor2}. We sample elements of $A$ from a uniform distribution over the unit interval from $0$ to $1$. In addition, we consider the following values for different parameters: F dimensionality $p$ is $100$, sample size $n$ changes from $10$ to $2000$, the sparsity level $k$ is 5, signal intensity $s$ ranges from $1$ to $5$. $w^*$ is picked at random from uniform distribution on $k$-sparse vectors and then scaled to have $\ell_1$ intensity of $s$. The recovery error is averaged for each $s$ and $n$ using $10$ Monte-Carlo Simulations. The averaged $\ell_2$ error is plotted against the upper bound given by Corollary \ref{Cor2} in Fig. \ref{Cor2Fig}. This figure demonstrates that the derived bound is almost tight as the points show almost a linear relationship between error and the upper bound.

\subsection{Rescaled LASSO vs. Regularized ML}
Parameter estimation based on LASSO for the Poisson setting has been studied in \cite{cite5}. The idea is to view the problem as an additive noise problem, where noise belongs to an exponential family of distributions. Alternatively, in \cite{cite5} the problem is viewed as an additive Gaussian noise problem with noise variance being equal to its mean to mimic ``Poisson like'' behavior. This results in a rescaled version of LASSO, which is then used to estimate model parameters. This amounts to scaling the loss function associated with each observation by the mean (or equivalently the variance). 

This approach motivates us to compare the regularized ML method against re-scaled LASSO for poisson distributed data, to highlight the essence of using regularized ML instead of rescaled LASSO for our setting. In this section we will demonstrate that our regularized ML outperforms re-scaled LASSO in several regimes including low SNR, high dimensions, and moderate to low sparsity levels.

To compare the performance of regularized ML and rescaled LASSO, we first generate a random sensing matrix $A \in \mathbb{R}^{n\times p}$ where each element $a_{i,j}$ is an independent truncated Gaussian random variable. According to Lemma \ref{RERa}, when the number of rows satisfies $n = o( k^2 \log( p ) \log^3 (k \log p ) )$, matrix $A$ satisfies RE with high probability. We also generate a random value  $\lambda_0$, such that $\lambda_{0,i}=\lambda_0$ for all $i$, and some sparse vector $w^*\in \mathbb{R}^p$, with $\|w^*\|_1=s$. To recover $w$, we generate $n$ Poisson distributed data with coefficients specified in $A$ as:
$$y_i=\Poisson(\lambda_{0}+a_i^{\T}w^*)$$
We first solve the non linear optimization where $w$ is constrained to be in $\Theta_s$. 
$$\widehat{w}_{ML}=\arg \min_{w \in \Theta_s}-\frac{1}{n}\sum_{i=1}^ny_i\log(\lambda_{0}+a_i^{\T}w)-a_i^{\T}w$$
For the purpose of comparison we compute the rescaled LASSO estimator.
$$\widehat{w}_{LS}=\arg\min_{w\in \Theta_s}\frac{1}{n}\sum_{i=1}^n\frac{(y_i-\lambda_{0}-a_i^{\T}w)^2}{\lambda_{0}+a_i^{\T}w}$$

\begin{figure}
\begin{minipage}[b]{1\linewidth}
  \centering
  \centerline{\includegraphics[width= 9cm]{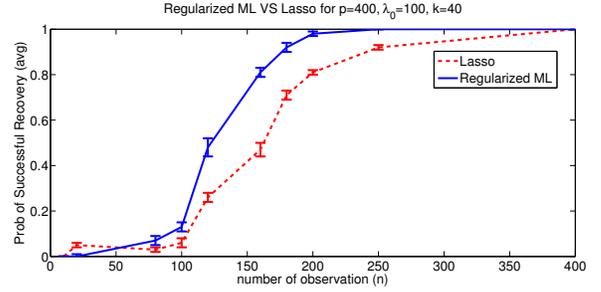}}
\caption{Probability of successful support recovery as a function of $n$ for $p=400$, $\lambda_0=100$, $k=40$, $t=10^{-4}$, and $m=100$ Monte Carlo loops. This figure illustrates twice faster convergence of probability of success with respect to number of observations for Regularized ML in comparison to Rescaled LASSO.
}
\label{tr}
\end{minipage}
\end{figure}

\begin{figure}
\begin{minipage}[b]{1\linewidth}
  \centering
  \centerline{\includegraphics[width= 9cm]{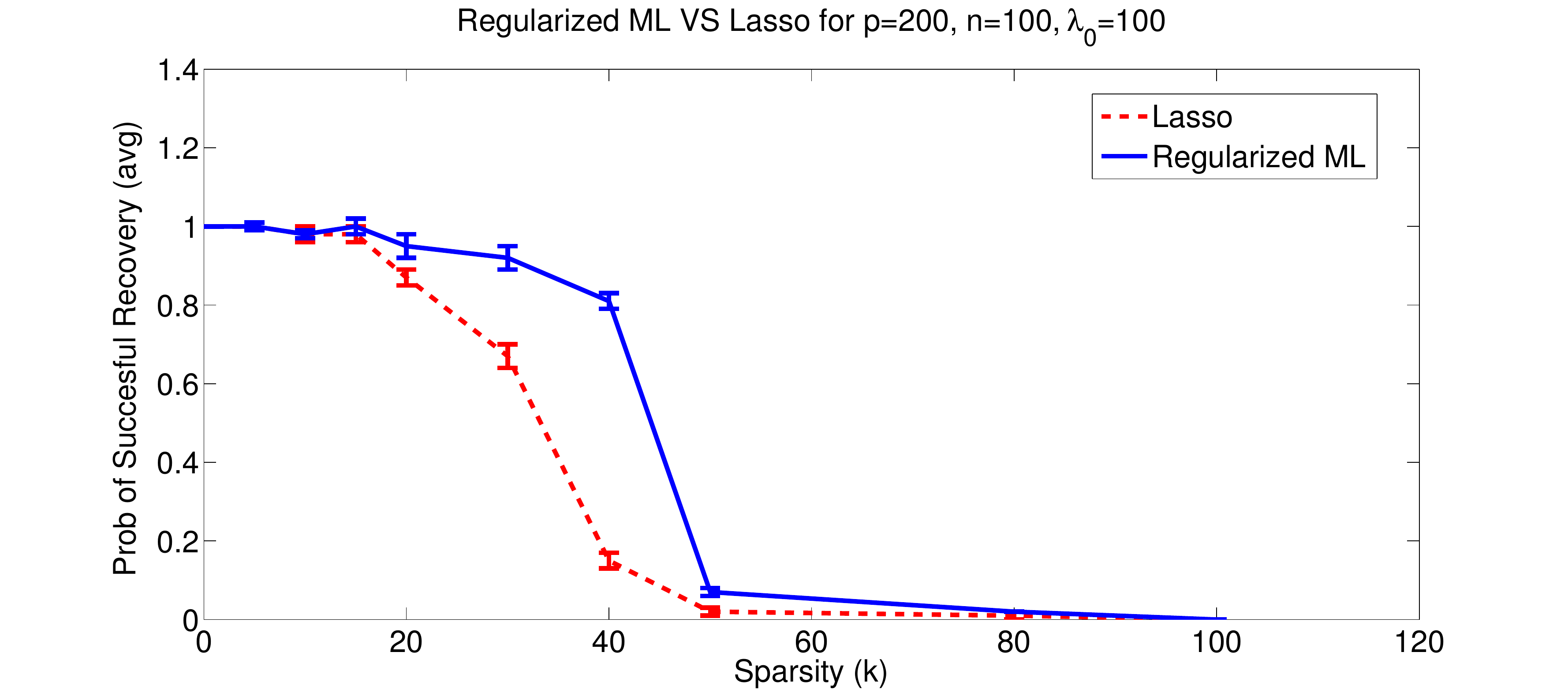}}
\caption{Probability of successful support recovery as a function of $k$ for $p=200$, $\lambda_0=100$, $n=100$, and $m=100$ Monte Carlo loops. As this figure suggests, probability of successfull recovery drops faster when the sparsity level increases for Rescaled LASSO in comparison to Regularized ML. This shows robustness of ML approach to model parameters.
}
\label{k}
\end{minipage}
\end{figure}

For comparison purposes we then threshold the solution by zeroing out components of $\widehat{w}_{ML}$ and $\widehat{w}_{LS}$ below a pre-defined small threshold $t$. 
We average the estimation performance over 100 Monte Carlo loops. The performance of the two methods are compared in Fig. \ref{tr} and Fig. \ref{k}. The results are compared in terms of number of observations $n$, and different sparsity levels $k$, respectively.

In Fig. \ref{tr}, we compare the result of regularized ML estimation with rescaled LASSO as a function of $n$. We fix $\lambda_0=100$, $p=400$, $t=10^{-4}$ and $k=40$. At each iteration, we estimate $w^*$ based on $n$ observations where $n$ varies from $2$ to $400$. We compare the performance of the two approaches based on probability of successful recovery of the support set. This error is 0 if the thresholded support set of the estimation is equal to that of the ground truth and 1 otherwise. We average this error over 100 samples of $w^*$ for a fixed $A$.

In Fig. \ref{k}, we compare the result of regularized ML estimation with rescaled LASSO for different sparsity levels, $k$. This time, we fix $\lambda_0=100$, $p=200$, and $n=100$. For each $k$, we generate 100 samples of $k$-sparse $w^*$'s and recover them from $n$ observations. Since $\|w^*\|_1=1$ for all values of $k$, we threshold each element of $w$ by $t = \frac{0.01}{k}$, to obtain their sparse support set. We measure the performance of the two estimators based on average probability of successful recovery of the thresholded support set for each value of $k$. 

Notice that the error bars in Fig. \ref{tr} and Fig. \ref{k} indicate that the difference between the methods is indeed statistically significant.

\begin{figure}
\begin{minipage}[b]{1\linewidth}
  \centering
  \centerline{\includegraphics[width= 9cm]{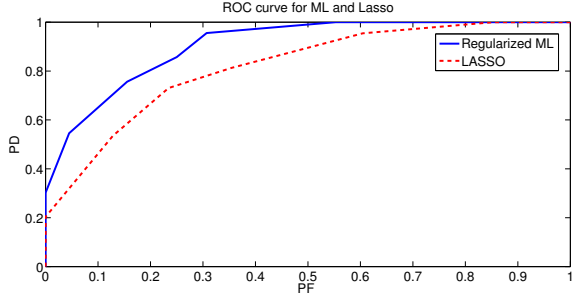}}
\caption{ROC curve for Regularized ML and Rescaled LASSO for  $n=100$, $p=200$, $k=20$, $\lambda=100$, and $m=100$ Monte Carlo loops. This figure illustrates the superiority of Regularized ML in comparison to Rescaled LASSO for parameter estimation under Poisson models.
}
\label{ROC}
\end{minipage}
\end{figure}

In Fig. \ref{ROC}, we compare the result of regularized ML estimation with rescaled LASSO in terms of the ROC curves. In an ROC curve, the average number of true detections is plotted against the average number of false alarms. True detections are indices that are common in the thresholded estimated support set and that of the Ground Truth, whereas, false alarms are the indices in the thresholded estimated support set that are not included in the support set of the Ground Truth. This time, we fix $\lambda_0=100$, $p=200$, $n=100$, and $k=20$. We fix a sensing matrix $A$, and generate 100 random $w^*$'s. By applying different thresholds $t=\frac{1}{k}$ to $t=\frac{0.001}{k}$ we obtain the different points in the ROC plot. We average Probability of Detection (PD) and Probability of False alarm (PF) over $100$ Monte Carlo loops.



\subsection{Explosive Identification}

In this experiment, we first measure the light intensities of different fluorophores before and after separate exposures to a unit weight of different explosives. The intensities are measured by counting the number of photons received at each photo-sensor. Each explosive $j$ has a unique quenching effect in the fluorescence property of each fluorophore $i$, which we denote by $a_{i,j}$. In the experimental setting, $\lambda_i$ is the before exposure intensity for fluorophore $i$ and is estimated by averaging the before exposure photon counts from multiple experiments. Therefore, the $\lambda_i$'s can be assumed to be known. We model the after exposure intensity $y_i$ as :
$y_i=\Poisson(\lambda_i(1-a_{ij}))$

In the next step, fluorophores are exposed to an unknown mixture of these explosives. The goal is to recover which and how much of each explosive is contained in that mixture.

The physics of the problem suggests that when the fluorophore is exposed to a mixture of explosives, the quenching effects are additive in the regime where the mixture weights are small \cite{cite6}. Therefore, our observations are best modeled by a Poisson distribution with additive rate model for each fluorophore:
$$ y_i \sim \Poisson\left(\lambda_i(1-\sum_{j=1}^pa_{ij}w^*_j)\right)$$ 
where $p$ is the total number of basic explosives and $w^*_j$ is the amount of the explosive $j$ in the mixture. We solve this problem through Regularized ML and Rescaled LASSO and compare the results.

\begin{figure}[htb]
\centering
\centerline{\includegraphics[width= 10cm]{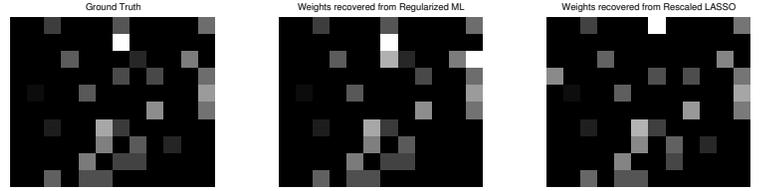}}
\caption{Sparse recovery results for $k\leq 3$. From left to right: Ground Truth, $w^*$, ML estimate of $w$ , LASSO estimate of $w$. Columns: Basic explosives.
Rows: Synthesized mixtures with $k$ basic explosives. Non-black squares at each row show the explosives that the corresponding mixture is composed of. lighter colors show larger amounts.
}
\label{FlComp}
\end{figure} 

In this problem, matrix $A$, the responses of $n=8$ fluorophores to $p=12$ basic explosives is given. Based on this given data and our additive model for mixtures, we generated 10 mixtures by combining up to 3 random explosives. We used Regularized ML and Rescaled LASSO to identify these mixtures through their effect on fluorescence property of our fluorophores. The result is shown in the form of a $10\times 12$ grid in Fig. \ref{FlComp}. In this grid, rows are different mixtures and columns are different explosives. Dark squares indicate the absence (or negligible contribution), whereas lighter squares indicate higher amount of the corresponding explosive in the associated mixture.

Since in this example photon count rates are of the order of $10^5$, Normal distribution could be considered as a good approximation to the corresponding Poisson distribution. Hence Regularized ML and Rescaled LASSO show similar behavior in this example.

\subsection{Internet Marketing Application}

In this application, our goal is to identify the most effective advertisement websites that result in higher website traffic in the clothing market. Our assumption is that the website traffic is generated as a superposition of the traffic generated from current customers and the traffic from advertisement through backward links (links in advertisement websites that are linked to these business websites). In general, big business websites typically buy a total of 1000-1500 backward links from a number of advertisement websites. However, the hypothesis is that only a few of these advertisement websites are efficiently directing costumers. Our goal is to identify those dominant advertisement websites.
 
We model the number of daily visits, $y_i$, by:
$$y_i=\Poisson(\lambda_{i,0}+a_i^{\T}w^*)$$
where $\lambda_{0,i}$ models the current customers who visit the site directly and is obtained through online statistics of the website. Specifically, a long run average of traffic which is not referred by the advertisement website gives a reasonable estimate and so we can assume that $\lambda_{0,i}$ to be known. This traffic could be logged and acquired through online statistics of the business websites. 

Moreover, $a_{i,j}$ is the number of backwards link for the website $i$ in the advertisement website $j$. Our model assumes that each of the backward links brings independent traffic to the website. Therefore, we used the Poisson distributed random variable described earlier to model the number of visits to a business website.

Our observations are the daily online visits to 50 top clothing brands. From the information provided in alexia.com, we chose the top 150 advertisement websites for these brands along with the number of backward links for each website. Our goal is to recover the weight vector $w^*$, where $w^*_j$ is a measure of dominance for advertisement website $j$ in clothing market. We recover $w^*$ via regularized ML and rescaled LASSO (weights smaller than 0.01 are theresholded to 0). The result is provided in Table 1. In this table, we illustrate the corresponding score for each popular website based on their dominance in advertising for clothing brands.

\begin{table}[h]
\caption{Top backwardlist websites for clothing brands using Regularized ML and Rescaled LASSO}
\begin{center}
\begin{tabular}{|l|c|c|c|c|c|c|c|c|c|}
               \hline
               Backward link & ML estimated weight \\
               \hline
               Amazon & 0.32\\
               Twitter & 0.21\\
               Pinterest & 0.17 \\              
               Google &  0.15\\
               Blogger &  0.06\\
               Bing & 0.05\\
               douban & 0.01\\
               tumblr & 0.01 \\              
               \hline
               \hline
               Backward link & Lasso estimated weight \\
               \hline
               Amazon & 0.35\\
               Pinterest & 0.17 \\   
               Twitter & 0.16   \\        
               Google &  0.16\\
               Bing & 0.13\\            
               \hline

\end{tabular}
\end{center}
\end{table}

To compare the result of the two approaches mentioned above, we use the Bayes factor \cite{Kass} and predictive held-out log likelihood comparison mentioned in \cite{Blei}. It should be mentioned that these tests are interpretable only when the number of parameters are comparable in the hypothesis models. In our problem, in fact, the two models have equal number of parameters.

Given a set of observed data $y_1, \hdots, y_n$, and a model selection problem in which we have to choose between two models, Bayesian inference compares the plausibility of the two different models $M_1$ and $M_2$ through a likelihood test:

$$\frac{\Pr\left(y_1, \hdots, y_n|M_1\right)}{\Pr\left(y_1, \hdots, y_n|M_2\right)}\lessgtr 1$$ 

When the parameters of models $M_1$ and $M_2$ are not known a priori, in Bayes factor test, we estimate them from $y_1, \hdots, y_n$ and then use those estimations in computing the likelihood ratio. On the other hand, in predictive held-out log likelihood comparison, we divide the data into two groups. We estimate the model parameters for $M_1$, and $M_2$ using the first group of data, and we compare the likelihoods for the second part of data given $M_1$ and $M_2$ specified by the first group.

Since Poisson is a PMF distribution on integers, to compare the two models using Bayes factor, we need to superimpose the Gaussian distribution on a histogram defined on integer valued $y_i$'s. For a Gaussian distribution characterized by $\mathcal{N}(\mu,\sigma)$, the value of the histogram at each integer valued $y$ is computed as:
\begin{equation}
\label{hist}
\hist(y)=\frac{1}{Q(\frac{\mu}{\sigma})}\times \left(Q\left(\frac{y-\mu}{\sigma}\right)-Q\left(\frac{y+1-\mu}{\sigma}\right)\right)
\end{equation}
It is easy to show that this histogram corresponds to a valid PMF. We denote this PMF by $\overline{\mathcal{N}}(\mu,\sigma)$. 

\begin{figure}
\begin{minipage}[b]{1\linewidth}
  \centering
  \centerline{\includegraphics[width= 9cm]{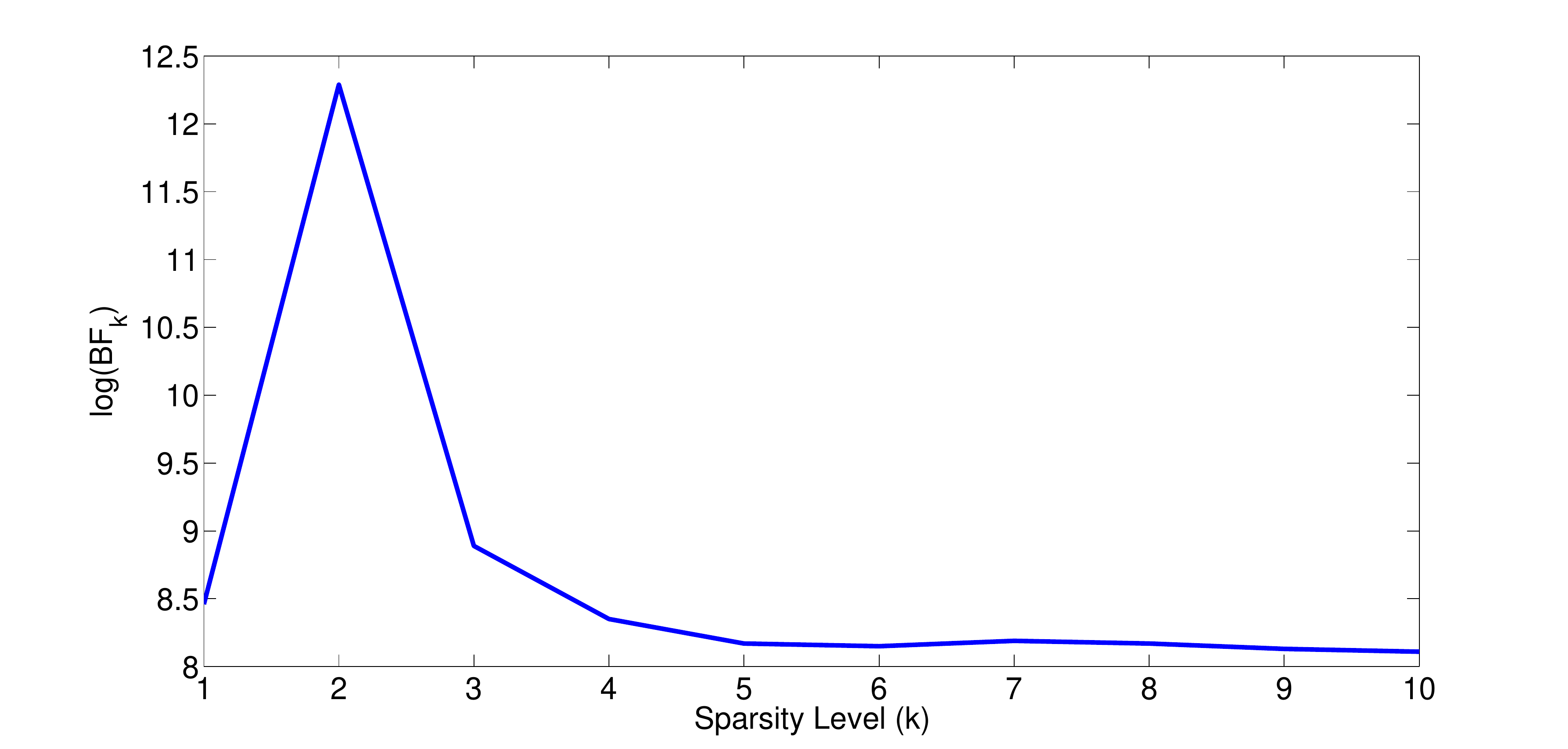}}
\caption{The Bayes Factor Ratio for regularized ML and rescaled LASSO. As we can see, higher Bayes Factor for Regularized ML suggests that Poisson is the right model for this problem.
}
\label{cox}
\end{minipage}
\end{figure}

\begin{figure}
\begin{minipage}[b]{1\linewidth}
  \centering
  \centerline{\includegraphics[width= 9cm]{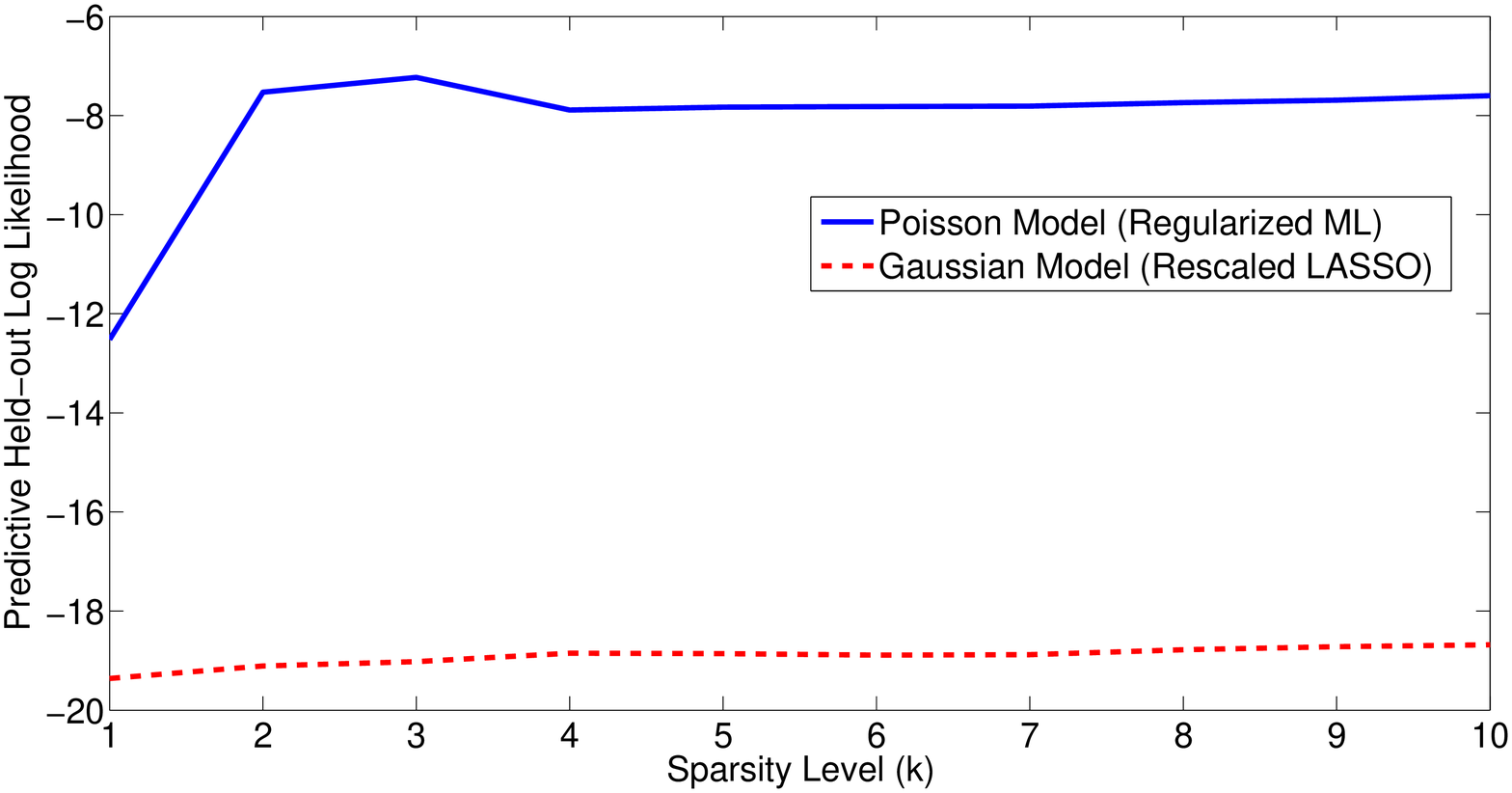}}
\caption{The Predictive Held-out Log Likelihood Comparison. The huge gap between the Predictive Held out Likelihood of Regularized ML and Rescaled Lasso implies that Poisson is the right model for this problem therefore, Regularized ML approach should be taken.
}
\label{cox1}
\end{minipage}
\end{figure}

After this conversion, the Bayes factor as a function of sparsity level, $k$, is calculated as:
$$BF_{k}=\frac{\Pr\left(y_1,\hdots, y_n|y_i\sim \Poisson(\lambda_0+A\widehat{w}^k_{ML})\right)}{\Pr\left(y_1, \hdots, y_n|y_i \sim\overline{\mathcal{N}}(\lambda_0+A\widehat{w}_{LS}^k,\lambda+A\widehat{w}^k_{LS})\right)}$$ 
where $\widehat{w}^k_{ML}$ and $\widehat{w}_{LS}^k$ are $k$ sparse theresholded approximations of $\widehat{w}_{ML}$ and  $\widehat{w}_{LS}$, respectively. The Bayes factor log curve as a function of sparsity is presented in Fig. \ref{cox}.
To compute the predictive held-out log likelihood for each method, we first use 80\% of the data (40 training data) to calculate an estimation of the parameters, $\widetilde{w}_{ML}$ and $\widetilde{w}_{LS}$. We use $\widetilde{w}_{ML}$ and $\widetilde{w}_{LS}$ for each model to compute the log likelihood function for the remaining 20\% of data (10 test data):
$$\mathcal{L}_{ML}=\sum_{i=1}^{10}-\lambda_i-a_i^{\T}\widetilde{w}_{ML}+y_i\log(\lambda_i+a_i^{\T}\widetilde{w}_{ML})-\log(y_i!)$$
\begin{align*}&\mathcal{L}_{LASSO}=\sum_{i=1}^{10}-\log\left(\!Q(\sqrt{\lambda_i+a_i^{\T}\widetilde{w}_{LS}})\right)+\\&\log\!\left(\!Q\left(\frac{y_i-\lambda_i-a_i^{\T}\widetilde{w}_{LS}}{\sqrt{\lambda_i+a_i^{\T}\widetilde{w}_{LS}}}\!\right)\!-\!Q\left(\frac{y_i+1-\lambda_i-a_i^{\T}\widetilde{w}_{LS}}{\sqrt{\lambda_i+a_i^{\T}\widetilde{w}_{LS}}}\right)\right)
\end{align*}
Intuitively, the model that is closer to the ground truth results in higher log likelihood value. The log likelihood values for the two approaches are shown in Fig. \ref{cox1}. The large gap between the predictive log likelihood of the two models implies that Poisson is a better underlying model for this application.
\subsection{Dynamics of Online Marketing}
In the previous section, our results show that ML estimator and Poisson model outperforms LASSO approach for the problem of online marketing. Therefore, in this section, we apply ML method to estimate the weights, $w^*$, for the advertisement websites over time.
\begin{figure}
\begin{minipage}[b]{1\linewidth}
  \centering
  \centerline{\includegraphics[width = 9cm]{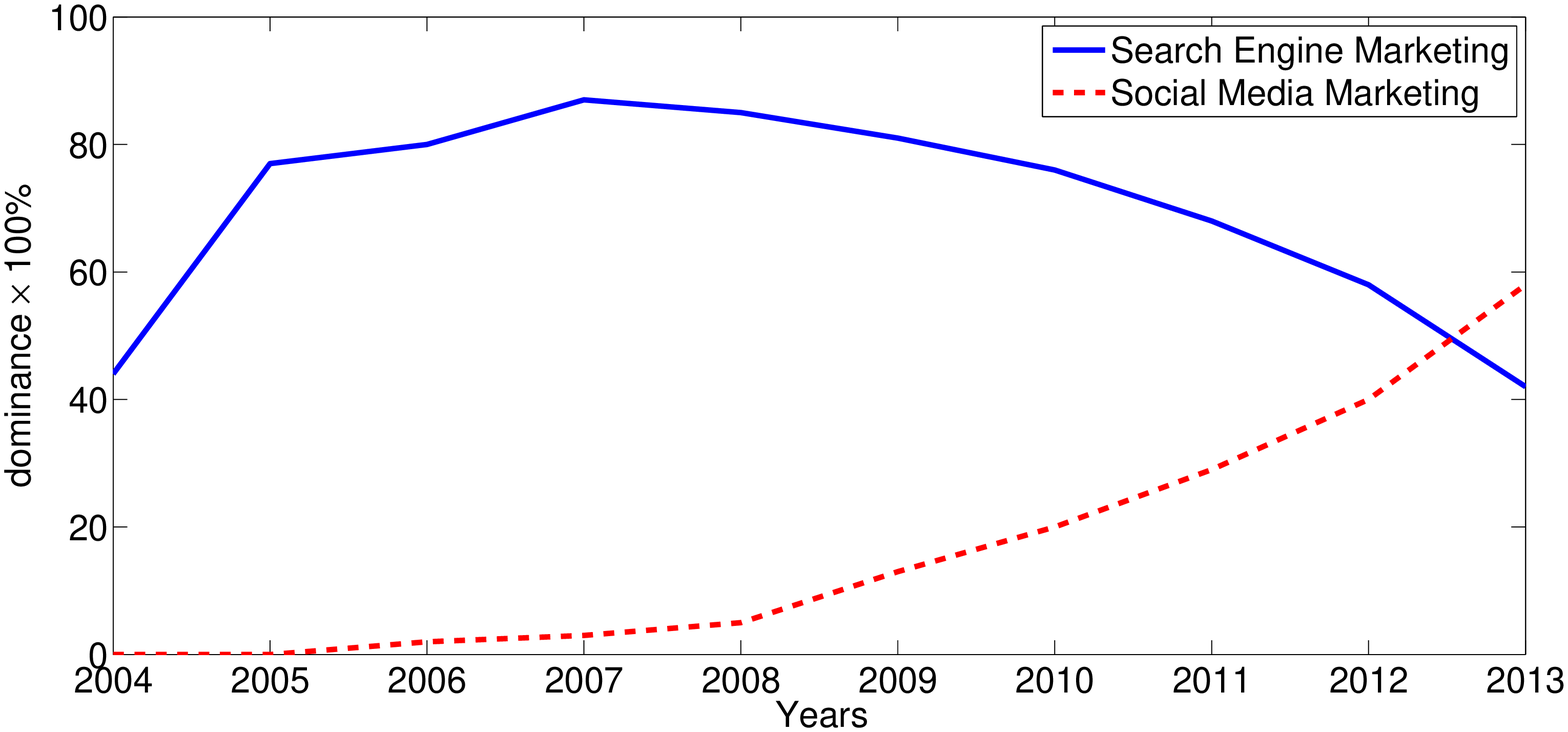}}
\caption{Dynamics of SEM and SMM over time for clothing market. This figure gives a quantitative comparison of the two most controversial methods of online marketing.
}
\label{cox2}
\end{minipage}
\end{figure}
\begin{table}[h]
\small{\caption{Top advertisement websites for clothing market in 2013 and 2008}
\begin{center}
\scalebox{0.7}{
\begin{tabular}{|l|c|c|c|c|c|c|c|c|c|}
               \hline
               Backward link & Estimated weights in 2013 &  Estimated weights in 2008\\
               \hline
               Twitter.com & 0.18 & 0.02\\
               Facebook.com & 0.18 & 0.02\\
               Pinterest.com & 0.14 & 0.00\\
               Amazon.com & 0.28 & 0.22\\              
               Google.com &  0.15 & 0.52\\
               Bing.com & 0.05 & 0.00\\
               Yahoo.com & 0.00 & 0.11\\
        
               \hline
\end{tabular}}
\end{center}}
\end{table}

A brief look at Table II shows how $w$'s have changed dramatically over time. To study this change closely, we estimated $w$'s for different advertisement websites from 2004 to 2013. We group the Social networks, such as facebook.com, twitter.com, pinterest.com, etc, together to study the effect of Social Media Marketing (SMM). We also group search engines, such as google.com, yahoo.com, bing.com, etc, together to represent Search Engine Marketing (SEM). We add the scores of the corresponding websites in each group. \hbox{Fig. \ref{cox2}} demonstrate the dynamics of SMM and SEM, the most controversial forms of online marketing, over time \cite{cite13}. Although SEM has been thought to be the most powerful media marketing tool, recent empirical studies show the growing influence of SMM during the last couple of years \cite{MF}.  The gigantic size of social media coupled with the relatively low cost per impression and the so called word of mouth have made SMM a powerful marketing tool. Our results confirm the significant influence of SMM relative to SEM since 2012.

\section{Conclusions}
We provided convergence guarantees for the solution of $\ell_1$ regularized ML decoder of a high dimensional sparse parameter for heterogeneous Poisson distributed data. 
Unlike least-squares linear regression setting, scale of the parameter has a significant effect on sample complexity. A new condition, Restricted Likelihood Perturbation (RLP), for successful recovery is introduced, which captures this effect. We then derived an expression relating RLP to the well-known restricted eigenvalue conditions. These expressions led us to deriving sample complexity bounds for several sensing matrix constructions. In our experiments, we verified the signal amplitude effect and tightness of the bounds. We also compared rescaled LASSO against our regularized ML to justify the essence of regularized ML in recovery based on our proposed observation model. 

\section{Appendix A} 
\label{AP}
\subsection{Useful Bounds:}
To show an exponential rate of convergence for Poisson distributed data, we need a tool to bound the tail probability. We build this tool from Bernstein inequality for Poisson distribution:
\begin{lemma}
\label{bert}
(Bernstein inequality)\cite{cite16} Let $y_1, \dots, y_n$ be independent random variables with means $\mu_1,\ldots, \mu_n$. Suppose that $\exists L>0$ such that $\forall k\in N$ and $k>1$:
$$\E[|(y_i-\mu_i)^k|]\leq \frac{1}{2}\E[(y_i-\mu_i)^2] L^{k-2}k! $$
Then, we have:
$$\Pr\left\lbrace\frac{1}{n}\sum_{i=1}^N |y_i-\mu_i|\geq \frac{2t}{n} \sqrt{\sum \E[(y_i-\mu_i)^2]}\right\rbrace \leq 2\exp{(-t^2)}$$

where $0<t\leq \frac{ \sqrt{\sum \E[(y_i-\mu_i)^2]}}{2L}$
\end{lemma}
\begin{lemma} 
\label{bert_cond}
For $y_i$'s distributed as:
$$y_i\sim \Poisson (\lambda_{i})$$
There exists a number $L > 0$, such that $\forall k\in N$ and $k>1$:
$$\E[|(y_i-\lambda_i)^k|]\leq \frac{1}{2}\E[(y_i-\lambda_i)^2] L^{k-2}k! $$ Moreover, $L = \max(1, \sqrt{\lambda_{max}})$.
\end{lemma}
\textbf{Remark: }The proof of this Lemma is partly provided in \cite{cite14} and is based on the fact that moment generating function for Poisson distribution with rate $\lambda$,
 $$M(t) = \exp\left(\lambda(\exp(t)-1) - \lambda t \right)$$
 is an analytic function, which means that its Taylor series converges around $t = 0$ on an open set in $\mathbb{R}$. Therefore, the $k$-th coefficient of Taylor series exists and is bounded. If $\lambda \leq 1$, it could be shown that all Taylor coefficients of $M(t)$ are less than or equal to $\frac{1}{2}$. If $\lambda > 1$, we replace $t$ by $t/\sqrt{\lambda}$. Then, all Taylor coefficients could be shown to be less than $\frac{1}{2}$. This completes the proof.

\subsection{\bf Proof of Lemma \ref{upperBnd}}
Now, we need to show that for any \hbox{$\epsilon>0$}, and the corresponding $\delta_{s,k}(\epsilon)$ defined in Eq. \eqref{deltask}, we have:
\begin{align*}
\Pr\left\{\lVert \widehat{w} - w^* \rVert_2\geq\epsilon\right\}
\leq   \Pr\left\{2\sup_{w \in \Theta_s}  | Q(w)-\overline{Q}(w) | \geq\delta_{s,k} \right\}
\end{align*}
The proof of this part, can be shown by combining two Lemmas:
\begin{lemma}
\label{lem2} For any \hbox{$\epsilon>0$}, we have:
\begin{align}
\label{bound2}
\Pr\left\{\lVert \widehat{w} - w^* \rVert_2 \geq \epsilon\right\}
\leq  \Pr\left\{| \overline{Q}(\widehat{w})-\overline{Q}(w^*) |  \geq \delta_{s,k} \right\}
\end{align}
\end{lemma}
\begin{lemma} 
\label{lem1}
For any $\delta>0$, we have:
\begin{align*}\Pr\left\{|\overline{Q}(\widehat{w})-\overline{Q}(w^*) |\!\geq\! \delta \right\}\leq\Pr\left\{\!\sup_{w \in \Theta_s}  | Q(w)-\overline{Q}(w) | \geq \frac{\delta}{2} \right\}
\end{align*}
\end{lemma}

{\it Proof of Lemma \ref{lem2}}:

Consider $$\text{Event } \mathcal{P} :  \lVert \widehat{w} - w^* \rVert_2\geq \epsilon$$
and
$$\text{Event } \mathcal{Q} :  | \overline{Q}(\widehat{w})-\overline{Q}(w^*) | \geq \delta_{s,k}(\epsilon) $$
We have:
$$(\mathcal{P} \Rightarrow \mathcal{Q})\Longrightarrow \Pr(\mathcal{P})\leq \Pr(\mathcal{Q})$$
Note that $\overline{Q}(w)$ is a convex function. Hence its minimum over a concave set $\lVert w - w^* \rVert_2 \geq \epsilon$ is on the boundary $\lVert w - w^* \rVert_2 = \epsilon$. Finally, by the definition in Eq. \eqref{deltask}, the minimum value on the boundary is lower bounded by $\delta_{s,k}$. Therefore, $\mathcal{P}$ implies $\mathcal{Q}$, which completes the proof. $\square$ \\

{\it Proof of Lemma \ref{lem1}}:

Based on the definition of $\overline{Q}(w)$:
$$w^*=\arg \min_{w \in \Theta_s}{\overline{Q}({w})}$$
Therefore:
\begin{align*}
|\overline{Q}(\widehat{w})-\overline{Q}(w^*)|
=  \overline{Q}(\widehat{w})-\overline{Q}(w^*)
\end{align*}
\begin{align*}
&\Pr\left\{|\overline{Q}(\widehat{w})-\overline{Q}(w^*)|\geq\delta\right\}\\ \nonumber
&=  \Pr\left\{\overline{Q}(\widehat{w})-{Q}(\widehat{w})+{Q}(\widehat{w})-\overline{Q}(w^*)\geq \delta\right\}\
\end{align*}
Note that according to the definition:
$${\widehat{w}}=\arg\min_{w\in\Theta_s}{Q}(w)$$
Therefore $Q({\widehat{w}}) \leq Q({{w}^*})$, and we have the following inequalities:
\begin{align*}
&\Pr\left\{\overline{Q}(\widehat{w})-{Q}(\widehat{w})+{Q}(\widehat{w})-\overline{Q}(w^*)\geq \delta\right\}\\\nonumber
&\leq  \Pr\left\{\overline{Q}(\widehat{w})-{Q}(\widehat{w})+{Q}({w}^*)-\overline{Q}(w^*)\geq \delta\right\}\\
&\leq\Pr\left\{2\sup_{w\in \Theta_s}|{Q}(w)-\overline{Q}(w)|\geq \delta\right\}
\end{align*}
which proves our claim.$\square$ \\

\subsection{\bf Proof of Lemma \ref{UltLem}}
We start with:
\begin{align}
\nonumber
&\Pr\left\{\sup_{w\in\Theta_s}|Q(w)-\overline{Q}(w)|\geq \frac{\delta}{2}\right\}\\
&\leq \Pr\left\{\sup_{w\in\Theta_s}\frac{1}{n}|\sum_{i=1}^n(y_i-\lambda_{w^*,i})\log(\lambda_{w,i})|\geq \frac{\delta}{2}\right\}\\
&\leq \Pr\left\{|\frac{1}{n}\sum_{i=1}^n(y_i-\lambda_{w^*,i})|\geq \frac{\delta}{2|\log(\lambda_{max})|}\right\}
\end{align}
To apply the result of Lemma $\ref{bert}$, we have to set:
\begin{equation}
\label{hadde_payin}
\frac{\delta}{2 |\log(\lambda_{max})|}= \frac{2t}{n} \sqrt{\sum_{i=1}^n\lambda_{w^*,i}}
\end{equation}
where 
\begin{equation}
\label{tt}
t\leq \frac{\sqrt{\sum_{i=1}^n\lambda_{w^*,i}}}{2L}
\end{equation} 
and $L := \max(1, \sqrt{\lambda_{max}})$.
Combining Eqn. $(\ref{hadde_payin})$ and Eqn. $(\ref{tt})$, we have: 
\begin{equation}
\frac{n\delta}{4|\log(\lambda_{max})|\sqrt{\sum_{i=1}^n\lambda_{w^*,i}}}\leq \frac{\sqrt{\sum_{i=1}^n\lambda_{w^*,i}}}{2L}
\end{equation}
Therefore, $\delta$ must be upper bounded by:
$$\delta\leq \frac{2\lambda_{min}|\log(\lambda_{max})|}{L}$$
to guarantee:
\begin{align}
\nonumber
&\Pr\left\{\sup_{w\in\Theta_s}|Q_n(w)-\overline{Q}_n(w)|\geq \frac{\delta}{2}\right\}\\
\label{balayi}
&\leq 2\exp\left(-\frac{n\delta^2}{16\log^2(\lambda_{max})\lambda_{max}}\right)
\end{align} 
Finally, note that as $\lambda_{\max} = \mathcal{O}(s)$, we may replace $\lambda_{max}$ by $s$ at the cost of getting some extra universal constants. 

\subsection{\bf Proof of Theorem \ref{mainTh2}}
Lets assume $u := \widehat{w} - w^*$. We have:
\begin{align}
\label{khkh}
\overline{Q}(\widehat{w})-\overline{Q}(w^*) = \underbrace{\frac{1}{n}\sum_{i=1}^n-\lambda_{w^*,i}\log\left(1+\frac{a_i^{\T}u}{\lambda_{w^*,i}}\right)+a_i^{\T}u}_{f(u)}
\end{align}
where $\lambda_{w^*,i}=\lambda_{0,i}+a_i^{\T}w^*$.

Our goal is to find a lower bound on $\delta_{s,k}$ defined earlier as:
$$\min_{\stackrel{\|u\|_2 = \epsilon}{\widehat{w}\in \Theta_s, ~~ w^{*} \in \Gamma_k \cap \Theta_s}} f(u) $$
We use the following inequality to lower bound $f$ :
\begin{equation}
a \leq b \Rightarrow a \log (1 + x/a) \leq b \log (1 + x/b)
\end{equation}
Hence
\begin{align}
\label{bound_f}
f(u)\geq &\frac{1}{n} \sum_{i=1}^n- \lambda_{max} \log\left(1+\frac{a_i^{\T}u}{\lambda_{max}}\right)+{a_i^{\T}u} \\
\geq & \frac{\lambda_{max}}{n} \sum_{i=1}^n- \log\left(1+\frac{a_i^{\T}u}{\lambda_{max}}\right)+\frac{a_i^{\T}u}{\lambda_{max}}
\end{align}
and from inequality $-\log(1+x)\geq -x$, we can show that for all $u$, $f(u)\geq 0$. 
Next, from Eqn. ($\ref{bound_f}$) we have:
\begin{align}
\label{chchch}
&\delta_{s,k} \geq \min\limits_{\stackrel{\|u\|_2=\epsilon}{\widehat{w} \in \Theta_s}, w^* \in \Gamma_k \cap \Theta_s} \frac{\lambda_{max}}{n} \sum_{i=1}^n-\log\left(1+\frac{a_i^{\T}u}{\lambda_{max}}\right)+\frac{a_i^{\T}u}{\lambda_{max}}
\end{align}
We make by the following change of variables:
$$X_i=\frac{a_i^{\T}u}{\lambda_{max}}$$
Before we proceed to apply the result of Definition 1, we need to check that $\|u_S\|_1\geq \|u_{S^c}\|_1$. We know:
$$\|u_S\|_1=\|\widehat{w}_S-w^*\|_1\geq \|w^*\|_1-\|\widehat{w}_S\|_1$$
Moreover, from $\|\widehat{w}\|_1\in \Theta_s$, we have:
$$\|\widehat{w}_S\|_1+\|\widehat{w}_{S^c}\|_1=\|\widehat{w}\|_1\leq s=\|w^*\|_1$$
Therefore,
$$\|u_{S^c}\|_1=\|\widehat{w}_{S^c}\|_1\leq \|w^*\|_1-\|\widehat{w}_S\|_1 \leq \|u_S\|_1$$
Now, from Definition 1, we have:
$$\frac{1}{n}\|X\|_2^2=\frac{1}{n}\sum_{i=1}^n X^2_i=\frac{1}{n}\sum_{i=1}^n \frac{(a_i^{\T}u)^2}{\lambda^2_{max}}\geq \gamma_k \frac{\epsilon^2}{\lambda^2_{max}} $$
Now, by applying Taylor series expansions around $X_i=0$ to each term in the sum in Eq. ($\ref{chchch}$), we have:
$$-\log\left(1+X_i\right)+X_i=-X_i+\frac{1}{(1+\widetilde{X}_i)^2}X^2_i+X_i$$
where $|\widetilde{X}_i|$ lies between 0 and $|X_i|$:
$$|\widetilde{X}_i|=|c\times 0+ (1-c)\times X_i|\leq |X_i|\leq \frac{2sa_{max}}{\lambda_{max}}$$ 
where the last inequality follows from the fact that $\|u\|_1\leq 2s$.

Therefore, we can rewrite Eqn. ($\ref{chchch}$) as
\begin{align*}
\delta_{s,k} &\geq\min_{\|X\|^2 \geq \gamma_k \frac{n\epsilon^2}{\lambda^2_{max}}}\frac{1}{n}\lambda_{max}\sum_{i=1}^n\frac{1}{(1+\frac{2sa_{max}}{\lambda_{max}})^2}X^2_i\\
&\geq \frac{\lambda_{max}}{(1+\frac{2sa_{max}}{\lambda_{max}})^2}\min_{\|X\|^2 \geq \gamma_k \frac{n\epsilon^2}{\lambda^2_{max}}} \frac{1}{n}\sum_{i=1}^n X^2_i\\
&\geq \gamma_k\frac{\lambda_{max}\epsilon^2}{\lambda^2_{max}(1+\frac{2sa_{max}}{\lambda_{max}})^2} \\ 
&= \gamma_k\frac{\lambda^3_{max}\epsilon^2}{\lambda^2_{max} (\lambda_{max} + 2 s a_{max} )^2} \\
&\geq \gamma_k\frac{\lambda^3_{max}\epsilon^2}{\lambda^2_{max} (\lambda_{max} + 2 \lambda_{max} )^2} \\
& = \frac{\gamma_k \epsilon^2}{9 \lambda_{max}}
\end{align*}

Finally, as $\lambda_{max} = \mathcal{O}(s)$, we replace $\lambda_{max}$ by $s$ and introduce some constants in the equation.

\bibliographystyle{IEEEbib}
\bibliography{refs}

\begin{thebibliography}{10}

\bibitem{cite6}
R.~C. Stringer, S.~Gangopadhyay, and S.~A. Grant,
\newblock ``Detection of nitroaromatic explosives using a fluorescent-labeled
  imprinted polymer,''
\newblock in {\em Analytical Chemistry}, 2010.

\bibitem{cite13}
J.~Beel, B.~Gipp, and E.~Wilde,
\newblock ``Academic search engine optimization (aseo): Optimizing scholarly
  literature for google scholar and co.,''
\newblock {\em Journal of Scholarly Publishing}, 2010.

\bibitem{cite20}
P.~J. Bickel, Ya, €™acov Ritov, and Alexandre~B. Tsybakov,
\newblock ``Simultaneous analysis of lasso and dantzig selector,''
\newblock {\em Annals of Statistics, Volume 37, Number 4 (2009), 1705-1732},
  2009.

\bibitem{Sara}
Sara A. Van~De Geer and Peter Buhlmann,
\newblock ``On the conditions used to prove oracle results for the lasso,''
\newblock {\em Electron. J. Stat}, 2009.

\bibitem{cite5}
J.~Jia, K.~Rohe, and B.~Yu,
\newblock ``The lasso under poisson-like heteroscedasticity,''
\newblock {\em arXiv}, 2010.

\bibitem{cite10}
S.~Negahban, P.~Ravikumar, M.~J. Wainwright, and B.~Yu,
\newblock ``A unified framework for high-dimensional analysis of m-estimators
  with desomposable regularizers,''
\newblock {\em arXiv}, 2012.

\bibitem{cite12}
I.~Rish and G.~Grabarnik,
\newblock ``Sparse signal recovery with exponential-family noise,''
\newblock {\em Allerton}, 2009.

\bibitem{cite14}
S.~Kakade, O.~Shamir, K.~Sridharan, and A.~Tewari,
\newblock ``Learning exponential families in high-dimensions: Strong convexity
  and sparsity,''
\newblock {\em arXiv.org}, 2009.

\bibitem{cite15}
S.~Portnoy,
\newblock ``Asymptotic behavior of likelihood methods for exponential families
  when the number of parameters tends to infinity,''
\newblock {\em Annals of Statistics}, 1988.

\bibitem{cite11}
M.~Raginsky, R.~M. Willett, Z.~T. Harmany, and R.~F. Marcia,
\newblock ``Compressed sensing performance bounds under poisson noise,''
\newblock {\em IEEE Transaction on Signal Processing, vol. 58, no. 8, pp.
  3990?4002}, 2010.

\bibitem{cite112}
M.~Raginsky R.~M.~Willett,
\newblock ``Compressed sensing performance bounds under poisson noise,''
\newblock {\em ISIT}, 2009.

\bibitem{cite113}
M.~Raginsky, S.~Jafarpour, Z.~T. Harmany, R.~F. Marcia, R.~M. Willett, , and
  R.~Calderbank,
\newblock ``Performance bounds for expander-based compressed sensing in poison
  noise,''
\newblock {\em IEEE Transaction on Signal Processing, vol. 59, no. 9, pp.
  4139?4153}, 2011.

\bibitem{Van09}
S.~van~de Geer and P.~Buhlmann,
\newblock ``On the conditions used to prove oracle results for the lasso,''
\newblock {\em Electron. J. Stat.}, vol. 3, pp. 1360--1392, 2009.

\bibitem{Rask10}
G.~Raskutti, M.~J. Wainwright, and B.~Yu,
\newblock ``Restricted eigenvalue properties for correlated gaussian designs,''
\newblock {\em Journal of Machine Learning Research}, vol. 11, pp. 2241--2259,
  2010.

\bibitem{ShuZh09}
Shuheng Zhou,
\newblock ``Restricted eigenvalue conditions on subgaussian random matrices,''
\newblock {\em arXiv:0912.4045v2 [math.ST]}, 2009.

\bibitem{Rud10}
M.~Rudelson and S.~Zhou,
\newblock ``Reconstruction from anisotropic random measurements,''
\newblock {\em arXiv:1106.1151 [math.ST]}, 2010.

\bibitem{cite17}
A.~Takeshi,
\newblock ``Advanced econometrics,''
\newblock {\em Harvard University Press}, 1985.

\bibitem{Kass}
R.~E. Kass and A.~E. Raftery,
\newblock ``Bayes factor,''
\newblock {\em Journal of the American Statistical Association, Vol. 90, No.
  430. (Jun., 1995), pp. 773-795}, 1995.

\bibitem{Blei}
D.~M. Blei and P.~I. Frazier,
\newblock ``Distance dependent chinese restaurant processes,''
\newblock {\em Journal of Machine Learning Research}, 2011.

\bibitem{MF}
M.~Wasiq A.~Bashar, I.~Ahmad,
\newblock ``Effectiveness of social media as a marketing tool: Anempirical
  study,''
\newblock {\em International Journal of Marketing, Financial Services and
  Management Research}, 2012.

\bibitem{cite16}
S.~Bernstein,
\newblock ``The theory of probabilities,''
\newblock {\em Gastehizdat Publishing House}, 1946.

\bibitem{cite1}
W.~K. Newey,
\newblock ``Uniform convergence in probability and stochastic equicontinuity,''
\newblock {\em Econometrica}, 1991.

\bibitem{cite8}
B.~Hoadley,
\newblock ``Asymptotic properties of maximum likelihood estimators for the
  independent not identically distributed case,''
\newblock {\em The Annals of Mathematical Statistics}, 1971.

\bibitem{cite2}
F.~Hayashi,
\newblock {\em Econometrics},
\newblock Princeton University Press, 2000.

\bibitem{cite3}
D.~L. Snyder,
\newblock {\em Random Point Processes},
\newblock Wiley-Interscience Publication, 1976.

\bibitem{cite18}
F.~X. Dupa, J.~Fadili, and J.~L. Starck,
\newblock ``Deconvolution under poisson noise using exact data fidelity and
  synthesis or analysis sparsity priors,''
\newblock {\em International Conference on Image Processing (ICIP)}, 2011.

\bibitem{cite19}
F.~X. Dupa, J.~Fadili, and J.~L. Starck,
\newblock ``A proximal iteration for deconvolving poisson noisy images using
  sparse representations,''
\newblock {\em IEEE Transaction on Image Processing}, 2011.

\end{thebibliography}
\newpage

\begin{IEEEbiographynophoto}{Delaram Motamedvaziri}
 received the B.Sc in Electrical Engineering from Sharif University of Technology, Tehran, Iran, in 2008. She is a PhD student in the Electrical and Computer Engineering Department at Boston
University. Her research interests are in High Dimensional Sparse Recovery, Statistical Signal Processing and Statistical Learning.
\end{IEEEbiographynophoto}

\begin{IEEEbiographynophoto}{Mohammad Hossein Rohban}
 received his Ph.D. in Computer Engineering (with focus on Artificial Intelligence) from Sharif University of Technology, Tehran, Iran,  in 2012. He is currently a postdoctoral associate in Electrical and Computer Engineering Department in Duke University.  His current research interests include high dimensional sparse recovery, and provable latent variable learning.
\end{IEEEbiographynophoto}

\begin{IEEEbiographynophoto}{Venkatesh Saligrama}
(SM’07) is a faculty member in the Electrical and Computer Engineering Department at Boston
University. He holds a PhD from MIT. His research interests are in Statistical Signal Processing, Statistical Learning, Video Analysis,
Information and Decision theory. He has edited a book on Networked Sensing, Information and Control. He has served as an Associate
Editor for IEEE Transactions on Signal Processing and Technical Program Committees of several IEEE conferences. He is the recipient of
numerous awards including the Presidential Early Career Award(PECASE), ONR Young Investigator Award, and the NSF Career Award.
More information about his work is available at http://blogs.bu.edu/srv
\end{IEEEbiographynophoto}
\IEEEpeerreviewmaketitle

\end{document}